\documentclass[11pt]{amsart}

\usepackage{mathrsfs}
\usepackage{amsmath, amscd, amsthm, amssymb, amsfonts, verbatim,subfigure}

\usepackage[all]{xy}
\usepackage[dvips]{graphicx,color}

\title[The Gromov-Witten potential associated to a TCFT]{The Gromov-Witten potential associated to a TCFT}
\author{Kevin Costello}
\address{Department of Mathematics \\University of Chicago}
\email{costello@math.uchicago.edu}

\date{12 September 2005}

\newcommand{\what}{\widehat}
\newcommand{\tr}{\triangle}

\newcommand{\K}{\mbb K}

\newcommand{\Comp}{\op{Comp}}

\newcommand{\br}{\overline}

\newcommand{\iso}{\cong}
\newcommand{\C}{\mathbb C}

\newcommand{\Q}{\mbb Q}

\newcommand{\norm}[1]{\left\| #1 \right\|}
\newcommand{\Oo}{\mathcal O}
\newcommand{\Z}{\mathbb Z}
\newcommand{\defeq}{\overset{\text{def}}{=}}
\newcommand{\into}{\hookrightarrow}

\newcommand{\Gr}{\operatorname {Gr}}
\newcommand{\op}{\operatorname}

\newcommand{\mbb}{\mathbb}
\newcommand{\mf}{\mathfrak}
\newcommand{\mc}{\mathcal}
\newcommand{\from}{\leftarrow}
\newcommand{\ip}[1]{\left\langle #1 \right\rangle}
\newcommand{\abs}[1]{\left| #1 \right|}
\newcommand{\cmod}{\overline {\mc M}}

\renewcommand{\d}{\mathrm{d}}
\renewcommand{\epsilon}{\varepsilon}

\newcommand{\Tate}{ {\op{Tate}}}

\newcommand{\Fk}{\mf F}
\newcommand{\dhat}{\what{\op d}}
\newcommand{\tilmod}{\widetilde{\mc M}}
\newcommand{\delbar}{ \br{\partial} }

\DeclareMathOperator{\Sym}{Sym}

\newtheorem*{utheorem}{Theorem}

\newtheorem{ntheorem}{Theorem}

\newtheorem*{conjecture}{Conjecture}

\newtheorem{theorem}{Theorem}[subsection]
\newtheorem{thm-def}{Theorem/Definition}[theorem]
\newtheorem{proposition}[theorem]{Proposition}

\newtheorem{lemma}[theorem]{Lemma}

\newtheorem{definition}[theorem]{Definition}

\numberwithin{equation}{subsection}

\newtheoremstyle{rem}
  {4pt}
  {5pt}
  {}
  {}
  {\itshape}
  {:}
  {3pt}
  {}

\theoremstyle{rem}

\newtheorem*{remark}{Remark}
\newtheorem*{remarks}{Remarks}

\begin{document}

\begin{abstract}
This is the sequel to my preprint``TCFTs and Calabi-Yau
categories''. Here we extend the results of that paper to construct,
for certain Calabi-Yau $A_\infty$ categories, something playing the
role of the Gromov-Witten potential.  This is a state in the Fock
space associated to periodic cyclic homology, which is a symplectic
vector space.  Applying this to a suitable $A_\infty$ version of the
derived category of sheaves on a Calabi-Yau yields the B model
potential, at all genera.

The construction doesn't go via the Deligne-Mumford spaces, but
instead uses the Batalin-Vilkovisky algebra constructed from the
uncompactified moduli spaces of curves by Sen and Zwiebach.  The
fundamental class of Deligne-Mumford space is replaced here by a
certain solution of the quantum master equation, essentially the
``string vertices'' of Zwiebach.  On the field theory side, the BV
operator has an interpretation as the quantised differential on the
Fock space for periodic cyclic chains.  Passing to homology,
something satisfying the master equation yields an element of the
Fock space.

\end{abstract}

\maketitle

\section{Notation}
We work throughout over a ground field $\K$ containing $\Q$. Often
we will use topological $\K$ vector spaces.  All tensor products
will be completed.  All the topological vector spaces we use are
inverse limits, so the completed tensor product is also an inverse
limit.

All the results remain true without any change if we work over a
differential graded ground ring $R$, and use flat $R$ modules. (An
$R$ module is flat if the functor of tensor product with it is
exact). We could also have only a $\Z/2$ grading on $R$.

\section{Acknowledgements}
I would like to thank Tom Coates, Ezra Getzler, Alexander Givental
and Paul Seidel for very helpful conversations, and Dennis Sullivan
for explaining to me his ideas on the Batalin-Vilkovisky formalism
and moduli spaces of curves.

\section{Topological conformal field theories}

Let $\mc S$ be the topological category whose objects are the
non-negative integers,  and whose morphism space $\mc S(n,m)$ is
the moduli space of Riemann surfaces with $n$ parameterised
incoming and $m$ parameterised outgoing boundaries, such that each
connected component has at least one incoming boundary.   These
surfaces are not necessarily connected.  Let $\mc S_{\chi}(n,m)
\subset \mc S(n,m)$ be the space of surfaces of Euler
characteristic $\chi$.
 $\mc S$ is a symmetric monoidal topological category,
under disjoint union.

Let $C_\ast$ be the functor of normalised\footnote{Normalised means
we quotient by the subcomplex of degenerate simplices.} singular
simplicial chains with coefficients in $\K$, any field containing
$\Q$.  As $C_\ast$ is a symmetric monoidal functor, $C_\ast(\mc S)$
is a differential graded symmetric monoidal category.  We also need
a shifted version : define $C_\ast^{(d)}(\mc S(n,m))$ by
$$
C_i^{(d)}(\mc S)(n,m) = \oplus_{\chi} C_{i + d(\chi - n + m)} (\mc
S_{\chi}(n,m))
$$
Note that $\chi - n + m$ is even, so this shift in degree doesn't
affect the signs.
\begin{definition}
A topological conformal field theory of dimension $d$ over $\K$ is a
symmetric monoidal functor $F : C_\ast^{(d)}(\mc S) \to
\op{Comp}_{\K}$ to the category of complexes over $\K$, with the
property that the tensor product maps $F(n) \otimes F(m) \to F(n+m)$
are quasi-isomorphisms.
\end{definition}
The notion of topological conformal field theory was introduced
independently by Getzler \cite{get1994} and Segal \cite{seg1999}.

One source of topological conformal field theories is the
following theorem.
\begin{utheorem}[C., \cite{cos_2004oc}]
Let $\mc C$ be a Calabi-Yau $A_\infty$ category of dimension $d$
over  $\mbb K$.  Then there is a topological conformal field
theory
 $F$, of dimension $d$,  with a natural
quasi-isomorphism $CC_{\ast-d}(\mc C)^{\otimes n}  \iso F(n)$,
where $CC_\ast$ refers to the Hochschild chain complex.
\end{utheorem}
\begin{remark}
To make the signs easier to deal with, I have changed the notation a
little from \cite{cos_2004oc}.  This explains the shift by $d$ in
the Hochschild chain complex, which wasn't present in
\cite{cos_2004oc}.
\end{remark}

We are also interested in $\Z/2$ graded TCFTs.  This is a symmetric
monoidal functor from $C_\ast(\mc S)$ to the category of $\Z/2$
graded complexes of vector spaces, compatible with differentials and
with the grading.   To keep notation simple, we will always work
with only a $\Z/2$ grading.

\section{Informal outline of the construction}

This paper constructs, for certain TCFTs, something playing the role
of the Gromov-Witten potential.  One way to do this is due to
Kontsevich \cite{kon2003}.  His idea is that, in certain
circumstances, we can extend the TCFT to include operations coming
from the Deligne-Mumford spaces.   Then the Gromov-Witten potential
is defined in the usual fashion, using the fundamental class of
Deligne-Mumford space and $\psi$ classes. However, it turns out that
there is a choice involved in this construction, essentially of a
trivialisation of the circle action on the TCFT. The present paper
provides an alternative to Kontsevich's construction, which is
canonical, but instead of a generating function gives us a state in
a Fock space  \footnote{ Previously, in an attempt to understand
Kontsevich's lecture \cite{kon2003}, I constructed a homotopy
functor which when applied to the uncompactified moduli spaces of
curves yields the Deligne-Mumford spaces.  It follows automatically
that applying the same functor to a TCFT yields something which
carries operations from the Deligne-Mumford spaces. I sketched this
construction in the introduction to \cite{cos_2004oc}, without
proof. However, one of the results stated there is not right.  My
calculation of the result of the functor applied to a TCFT was
over-optimistic.  I had claimed we get cyclic homology, but instead
we get cyclic homology tensored with the ring of functions on a
certain space of inner products on cyclic homology. So that in order
to get operations from Deligne-Mumford spaces we need to choose such
an inner product. }.


 The constructions in this paper
bypass Deligne-Mumford space completely. Instead we use Sen and
Zwiebach's \cite{sen_zwi1994,sen_zwi1996} Batalin-Vilkovisky
algebra, which is constructed from the uncompactified moduli spaces
of curves.  The fundamental class of Deligne-Mumford space is
replaced by a certain solution of the quantum master equation in
this BV algebra. We construct from this solution of the master
equation a ray in a certain Fock space associated to a TCFT.

The idea that the fundamental chain satisfies the master equation is
not new here, but is due to Sullivan \cite{sul2004}, and appeared
implicitly in earlier work of Sen and Zwiebach \cite{sen_zwi1994,
sen_zwi1996}. This fundamental chain is essentially what Zwiebach
calls ``string vertices''; it is unique up to homotopy.

The connection with the Fock space seems to be new, though. There is
also a BV algebra associated to a TCFT.  The solution of the master
equation in moduli spaces gives one, say $S$,  in this BV algebra.
The master equation says that
$$
(\d + \tr) \exp S = 0
$$
We interpret the total BV operator $\d + \tr$ as the quantised
differential on a  chain level Fock space for a certain dg
symplectic vector space.  With this differential, the Fock space
becomes a dg module for the Weyl algebra. Thus, passing to homology,
$\exp S$ becomes an element of the Fock space for the homology of
our symplectic vector space.

\vspace{5pt}

There are various technicalities which make a rigourous exposition
of this construction a little unreadable. Therefore I'll start by
giving a sketch of the construction, which emphasises the geometry
and de-emphasises the technical details.

\subsection{Complexes with a circle action}
Let $F$ be a TCFT.  In this section we will make the simplifying
assumption that the maps $F(1)^{\otimes n} \to F(n)$ are
isomorphisms, and not just quasi-isomorphisms.  This is just for
expository purposes.

Let $$V = F(1)$$  Then $F(n) = V^{\otimes n}$.

As the monoid $\mc S(1,1)$ contains $S^1$ as a subgroup\footnote{We
allow ``infitely thin'' annuli in $\mc S(1,1)$, so that $\mc S$
becomes a unital category.}, the algebra $C_\ast(S^1)$ acts on $V$.
This is formal; there is a quasi-isomorphism \footnote{This only
works if $C_\ast$ is the normalised singular simplicial chains.}
$H_\ast(S^1) \to C_\ast(S^1)$  .   Let $\op D : V \to V$ be the odd
operator corresponding to the fundamental class of $S^1$.   We have
three associated complexes,
\begin{align*}
V_{\Tate} &= V((t)) \\
V^{h S^1} &= V[[t]] \\
V_{hS^1} &= V_\Tate / V^{hS^1} = t^{-1} V[t^{-1}]
\end{align*}
each with differential
$$
\d ( v \otimes f(t) ) = \d v \otimes f(t) + \op D v \otimes t f(t)
$$
The subscript or superscript $hS^1$ refers to homotopy coinvariants
on invariants respectively.
\subsection{The category of annuli}
Let $\mc M(m)$ be the moduli space of Riemann surfaces with $m$
(outgoing) boundary components.  Such surfaces may be disconnected;
also they may have connected components with no boundary. We allow
$m = 0$. Let $\mc M_g(m) \subset \mc M(m)$ be the subspace of
connected surfaces of genus $g$.

Now we define a topological category $A$, which is a subcategory of
$\mc S$.  The objects of $A$ are the non-negative integers, and the
morphisms are the morphisms in $\mc S(n,m)$ given by Riemann
surfaces all of whose connected components are annuli.  As each such
annulus has at least one incoming boundary component, this is a
subcategory.  Sending $m \mapsto \mc M(m)$ defines a symmetric
monoidal functor $A \to \op{Top}$. \footnote{By taking the coarse
moduli space of the orbispace $\mc M(m)$. } The maps $A(m,n) \times
\mc M(m) \to \mc M(n)$ are given by gluing annuli onto the boundary
of the surfaces in $\mc M(m)$.

Taking singular chains, we find a functor $C_\ast(A) \to
\op{Comp}_\K$, sending $n \mapsto C_\ast(\mc M(n))$.

If $F$ is a TCFT, sending $n \mapsto F(n)$ also defines a functor
from $C_\ast(A)$ to complexes.  A TCFT contains operations from the
uncompactified moduli spaces of curves, so we should be able to
relate $F$ and $\mc M$.  One could hope for a natural transformation
$C_\ast(\mc M) \to F$, encoding the TCFT structure. However, there
is a problem; $F$ carries operations from the moduli space $\mc
S(n,m)$ of Riemann surfaces which have at least one incoming
boundary component, whereas $\mc M$ is given by Riemann surfaces
with no incoming boundary components, and possibly no boundary at
all.

The  way around this problem is to construct a kind of semi-direct
product functor $F^{\mc M} : C_\ast(A) \to \Comp_\K$, which contains
both $F$ and $\mc M$, as well as the data of the action of $\mc S$
on $F$.  This construction will be explained in the body of the
paper. Here we will simply pretend that there is a natural
transformation $C_\ast(\mc M) \to F$.  This is purely for expository
purposes.  The argument used in the body of the text is more
complicated but relies on the same ideas.

\subsection{Batalin-Vilkovisky algebras}

Recall \cite{get1994} that a Batalin-Vilkovisky algebra is a
commutative dga $B$ equipped with an odd differential operator
$\tr$, which is of order $2$, and satisfies $\tr^2 = [\op{d},\tr] =
0$.

Consider the complex
$$
\Fk(\mc M) = \oplus C_\ast \left(\mc M(m) / S^1 \wr S_m \right)
$$
where $S^1 \wr S_m$ is the wreath product group $(S^1)^m \ltimes
S_m$.  Let $\Fk^{g,n}(\mc M) $ be the part coming from connected
surfaces of genus $g$ with $n$ boundaries.

The complex $\Fk(\mc M)$ is a  commutative differential graded
algebra, where the product comes from disjoint union of surfaces. We
will give $\Fk(\mc M)$ the structure of Batalin-Vilkovisky algebra.
This BV algebra from moduli spaces was introduced by Sen and
Zwiebach \cite{sen_zwi1994, sen_zwi1996}, and also studied by
Sullivan \cite{sul2004}.

The operator
$$\tr : \Fk_i(\mc M) \to \Fk_{i+1}(\mc M)$$
is defined\footnote{This operator is apparently not canonically
defined, as pull-back is not well defined on simplicial chains.
However, I'm lying about various technical issues.  We really should
be using $C_\ast(\mc M(m))_{h S^1 \wr S_m }$, which is take the
homotopy coinvariants by the action of the algebra $H_\ast(S^1) \wr
S_m)$, and then the operator is well defined. } as a sum of all
possible ways of gluing pairs of boundary components together, with
a full twist by $S^1$.
 This operator comes from the fundamental one chain
in $C_1 (A(2,0))$.

It turns out that any other symmetric monoidal functor $C_\ast(A)
\to \Comp_\K$ defines a Batalin-Vilkovisky algebra, in a similar
way. In particular, there is a BV structure on $\Sym^\ast V_{hS^1}$.
A point in $A(2,0)$, thought of as a zero chain, gives a pairing
$\ip{\ , \ }$ on $V$.  The operator $\tr$ on $\Sym^\ast V_{hS^1}$ is
the order $2$ differential operator, which on $\Sym^{\le 1}
V_{hS^1}$ is zero, and satisfies
$$
\tr((v_1 f_1(t_1))(v_2 f_2(t_2)) ) = \ip{\op D v_1, v_2 } \op{Res}
f_1 f_2 \op d t_1 \op d t_2
$$
Here we identify $V_{hS^1}$ with $V \otimes t^{-1} \K [t^{-1}]$, and
$v_i \in V$, $f_i \in t^{-1} \K[t^{-1}]$.

There is a map of BV algebras
$$
\Fk(\mc M) \to \Sym^\ast V_{hS^1}
$$
arising from the natural transformation $C_\ast(\mc M) \to F$.

If $B$ is any BV algebra, an element $S \in B$ satisfies the quantum
master equation if $(\d + \tr)\exp(S) = 0$.

This definition extends to elements $S \in B \otimes R$ for any ring
$R$.  We are interested in series $S \in B[[\lambda]] = B \otimes
\K[[\lambda]]$.
\begin{ntheorem}
There exists a sequence of elements $S_{g,n} \in \Fk^{g,n} (\mc M
)$, of degree $6g-6+2n$, with the following properties.
\begin{enumerate}
\item
 $S_{0,3}$ is a  $0$-chain of degree $1/3!$ in the moduli
space of Riemann surfaces with $3$ unparameterised, unordered
boundaries.
\item
Form the generating function
$$
S = \sum_{\substack{ g,n \ge 0 \\ 2g-2+n > 0}} S_{g,n}
\lambda^{2g-2+n}  \in  \lambda \Fk(\mc M) [[\lambda]]
$$
$S$ satisfies the Batalin-Vilkovisky quantum master equation :
$$
(\op{d} +  \tr) e^{S } = 0
$$
\end{enumerate}
Further, such an $S$ is unique up to homotopy through such series.
In particular the class $[e^{S}]$ in $\d + \tr$ homology is
independent of any choice.
\end{ntheorem}
A homotopy of solutions of the master equation (or of anything else)
is a family of such, parameterised by the contractible dga
$\K[t,\op{d}t] = \Omega^\ast_{\mbb A^1}$.

\begin{remarks}
\begin{enumerate}
\item
This result is essentially a mathematical formalisation of the work
of Sen and Zwiebach \cite{sen_zwi1994, sen_zwi1996}.  The choice of
such a solution of the master equation is essentially  the same as
the choice of string vertices in their work, which for these authors
is a certain subspace of $\mc M_g(n)$.  They realised that string
vertices satisfy the master equation, and that changing the choice
of string vertices changes $e^S$ by a $\d + \tr$ exact term. Also
the string vertices must correspond to the fundamental class, since
every Riemann surface appears uniquely by gluing surfaces lying in
the string vertices. This point was made clear in the work of
Sullivan \cite{sul2004} on chain level Gromov-Witten invariants.
\item
The proof of this theorem is very easy.  The only facts we need
about moduli spaces are trivial facts about the rational homological
dimension of $\mc M_{g,n}/S_n$.    Therefore the result is true if
we use any other sequence of spaces $\mc M_g(n)$ with the same
gluing structure.   The solution of the master equation is intrinsic
to the modular operad with compatible circle actions given by the
spaces $\mc M_g(n)$.
\end{enumerate}
\end{remarks}

Let $\cmod(n)$ be the moduli space of stable, possibly nodal,
possibly disconnected, algebraic curves, with $n$ marked smooth
points.  Consider the commutative dga $ \Fk(\cmod)$ defined by
$$
\Fk (\cmod) = \oplus_n C_\ast(\cmod(n)/S_n)
$$
The algebra structure comes from disjoint union.

Make this into a BV algebra by setting $\tr = 0$.  Let
$$
[\cmod] = \sum_{g,n} [\cmod_{g,n}/S_n] \lambda^{2g-2+n} \in
\Fk(\cmod)[[\lambda]]
$$
where $[\cmod_{g,n}/S_n]$ is an orbifold fundamental chain for
$\cmod_{g,n}/S_n$.
\begin{ntheorem}
There is a map $\Fk(\mc M)\to \Fk(\cmod)$ in the homotopy category
of BV algebras, such that $S$ maps to $[\cmod]$.
\end{ntheorem}

These results are not as mysterious as they might at first appear.
The master equation can be rephrased as
$$
\op{d} S_{g,n} +  \sum_{\substack {g_1+g_2 = g\\n_1+n_2=n+2}}
\frac{1}{2} \{ S_{g_1,n_1} , S_{g_2,n_2} \} + \tr S_{g-1,n+2} = 0
$$
where $\{ \quad \}$ is a certain odd Poisson bracket on the space
$\Fk(\mc M)$, constructed in a standard way from the BV operator
$\tr$.

If $\alpha \in \Fk^{g,n}(\mc M)$ and $\beta \in \Fk^{h,m}(\mc M)$
then $\{\alpha,\beta\}\in\Fk^{g + h,n+m-2}(\mc M)$ is the sum over
ways of gluing a boundary of $\alpha$ to one of $\beta$, with a full
twist by $S^1$. The twist by $S^1$ raises degree by $1$.  Similarly,
$\tr \alpha$ is the sum over ways of gluing a boundary of $\alpha$
to itself, with a twist by $S^1$.

We can see why the fundamental chain satisfies the master equation,
and relate it to the fundamental class of Deligne-Mumford space,
using a nice model for the spaces $\mc M(m)$ introduced by Kimura,
Stasheff and Voronov \cite{kim_sta_vor1995}. This is the moduli
space $\mc N(m)$ of algebraic curves in $\cmod(m)$ decorated at each
marked point with a ray in the tangent space, and at each node with
a ray in the tensor product of the rays of the tangent spaces at
each side. This is an orbifold with corners, whose interior consists
of non-singular curves; this shows that it is homotopy equivalent to
$\mc M(m)$. The operations of gluing marked points and of rotating
the rays at the marked points make $\mc N(n)$ into a functor from a
category $A'$, weakly equivalent to $A$.

We can construct from $\mc N$ a BV algebra $\Fk(\mc N)$, in the same
way we constructed $\Fk(\mc M)$.  These two BV algebras are
quasi-isomorphic.

Consider the space
$$X(n) \defeq   \mc N(n) /   S^1 \wr S_n $$
where $S^1 \wr S_n$ refers to the wreath product group $(S^1)^n
\ltimes S_n$.

A surface in $X(n)$  has unordered marked points, and
unparameterised boundaries (i.e. no ray in the tangent space).

The BV algebra $\Fk(\mc N)$ is given by
$$
\Fk(\mc N) = \oplus_n C_\ast(X(n))
$$
The space $X(n)$ is an orbifold with corners. Let $X_g(n)$ be the
subspace of connected surfaces of genus $g$.  The boundary of
$X_g(n)$ is (away from codimension $2$ strata) a union of bundles
over products of similar moduli spaces. There is a component for
each way of splitting $g = g_1+g_2$, $n +2 = n_1 + n_2$, and a
component corresponding to the loop, where we have a genus $g-1$
surface with $n+2$ marked points. Let us describe this last
component in detail.  It is a bundle over the space $X_{g-1}(n+2)$,
consisting of a point in this moduli space together with a choice of
two marked points, and a way of gluing them together. (There is an
$S^1$ of possible ways of gluing). Let us call this space $Y$. There
is a diagram
$$
X_{g-1}(n+2)  \from Y \to X_g(n)
$$
Pulling a chain in $X_{g-1}(n+2)$ back to $X$ (and so increasing the
degree by $1$),  and then pushing it forward to $X_g(n)$, is
precisely the operator $-\tr$.

A similar picture holds at the other boundary components, except we
find the bracket operator $\{\quad\}$ instead of $\tr$.

Now suppose that the orbifold with corners $X_g(n)$ has a
fundamental chain $[ X_g(n) ]$, behaving in a nice functorial
manner.  Then the discussion above would imply that
$$
\op{d} [X_g(n)] +  \sum_{\substack {g_1+g_2 = g\\n_1+n_2=n+2}}
\frac{1}{2} \{ [X_{g_1}(n_1)] ,[ X_{g_2(n_2)}] \} + \tr
[X_{g-1}(n+2)] = 0
$$
In other words, the fundamental chains of these moduli spaces
satisfy the BV master equation.

Now it is not difficult to show (using the fact that $H_i ( X_{g}(n)
) = 0$ for $i \ge 6g-7+2n$, $(g,n) \neq (0,3)$) that there is a
unique solution of the master equation up to homotopy, which sits in
the correct degrees and has the correct leading term.  Also,
quasi-isomorphic BV algebras have the same set of homotopy classes
of solutions of the master equation. These results now explains why
the image of the class $S_{g,n}$ in $H_{6g-6+2n}(\cmod_{g,n}/S_n)$
is the fundamental class.

\bigskip

\subsection{Weyl algebras and  Fock spaces}

We have a map of BV algebras  $\Fk(\mc M) \to \Sym^\ast V_{hS^1}$.
Let $\mc  D\in \Sym^\ast V_{hS^1} [[\lambda]]$ denote the image of
$\exp S$. As $S$ satisfies the master equation, $\mc D$ satisfies
$$(\d + \tr) \mc D = 0$$
The last step involves interpreting the homology class of $\mc D$ as
an element of a Fock space, by interpreting the differential $\d +
\tr$ as the natural differential on a chain-level Fock space.

$V$ carries an inner product, $\ip{\quad}$, coming from an annulus
in $C_\ast (\mc S(2,0))$.  We can arrange so that $\op D$ is skew
self adjoint, and $\d $ is self adjoint.

Define an antisymmetric pairing $\Omega$ on $V_{\Tate}$ by
$$
\Omega ( v f(t), w g(t) ) = \ip{v,w} \op {Res} f(-t) g(t) \op d t
$$
$\Omega$ is compatible with the differential $\d + t \op D$ on
$V_\Tate$.  This is the same form as that used in the work of Coates
and Givental \cite{coa_giv2001,giv2001,giv2004}.  The symplectic
nature of Tate cohomology is also studied by Morava \cite{mor2001}.

Let $\mc W(V_\Tate)$, the Weyl algebra, be the free algebra
generated by $u \in V_\Tate$ modulo the relation
$$
[u,u'] = \Omega(u,u')
$$
Here $[u,u']$ is the super commutator.

$V_\Tate$ has a decomposition $V_\Tate = V^{hS^1} \oplus V_{hS^1}$
into isotropic subspaces.  The subspace $V^{hS^1}$ is preserved by
the differential, but $V_{hS^1}$ is not in general.  The left ideal
in the Weyl algebra generated by $V^{hS^1}$ is also preserved by the
differential. Let $\Fk(V_\Tate)$ be the quotient module. The action
of $\Sym^\ast V_{hS^1} \subset \mc W(V_\Tate)$ on the image of $1$
in $\Fk(V_\Tate)$ gives an identification
$$
\mf F(V_\Tate) = {\Sym}^\ast V_{hS^1}
$$
The Weyl algebra action here is such that $V_{hS^1}$ acts by
multiplication, and $V^{hS^1}$ acts by differentiation.

$\Fk(V_\Tate)$ is a dg module for the dg algebra $\mc W(V_\Tate)$,
i.e.\ the differential is compatible with the action.
\begin{lemma}
The differential on $\Fk(V_\Tate)$ is $\d + \tr$. \label{lemma
quantised bv}
\end{lemma}
\begin{proof}
We can consider $\Sym^\ast V_{hS^1} \subset \mc W(V_\Tate)$ as a
subalgebra, which is not preserved by the differential.  The
differential $\op d_{\mc W(V_\Tate)}$ is a derivation.  So the
associated map $\op d_{\mc W(V_\Tate)} : \Sym^\ast V_{hS^1} \to \mc
W(V_\Tate)$ is characterised by how it behaves on the generators
$V_{hS^1}$. We have
$$
\op d_{\mc W(V_\Tate)} (v f(t)) = \op d_{V_{hS^1}} ( v  f(t) ) + \op
D v \op{Res} f(t) \op d t
$$
for $v f(t) \in V_{hS^1}$.  Note that for all $w g(t) \in V_{hS^1}$,
$$
\Omega ( \op {D} v \op{Res} f(t) \op d t, w g(t)    ) = \Omega (
\op{D} v tf(t) , w g(t) )
$$
Therefore, using the relations in $\mc W(V_\Tate)$, we find
\begin{multline*}
\op d_{\mc W(V_\Tate)} (v_1 f_1 \cdot v_2 f_2 \cdots v_k f_k ) = \op
d_{\Sym^\ast V^{hS^1}} ( v_1 f_1 \cdot v_2 f_2 \cdots v_k f_k ) \\ +
\sum_{i < j} v_1 f_1 \cdots \what{v_i f_i }  \cdots \what{ v_j f_j}
\cdots v_k f_k \Omega (\op D v_i t f_i, v_j f_j   ) \\ + \text{
terms in the left ideal generated by } V^{hS^1}
\end{multline*}
where $\what{ \ }$ indicates that we skip that term.

Modulo the ideal generated by $V^{hS^1}$ this is the same as $\d +
\tr$.
\end{proof}

The operator $\d + \tr$ is the unique odd differential on the space
$\Sym^\ast V_{hS^1}$ which make it into a dg module for the Weyl
algebra.

Let us use the notation
\begin{align*}
\mc H &= H_\ast(V_\Tate) \\
\mc H_+ &= H_\ast(V^{hS^1})
\end{align*}
We will assume the map $\op D : V \to V$ is zero on
homology\footnote{When the TCFT comes from a Calabi-Yau $A_\infty$
category, this corresponds to degeneration of the non-abelian Hodge
to de Rham spectral sequence}, and that the pairing on $H_\ast(V)$
is non-degenerate.   This implies that $\mc H$ is symplectic, and
the map $\mc H_+ \to \mc H$ is injective with Lagrangian image.

It is easy to see that
\begin{align*}
H_\ast (\mc W( V_\Tate) ) &= \mc W (\mc H) \\
H_\ast ( \Fk ( V_\Tate) ) &= \Fk ( \mc H)
\end{align*}
where $\Fk(\mc H)$ is the quotient of $\mc W(\mc H)$ by the left
ideal generated by $\mc H_+$.

Therefore
$$
[\mc D] \in \Fk ( \mc H)
$$
is an element in the Fock space for $\mc H$.

Of course, this construction is not restricted to the fundamental
class.   Suppose $R$ is a graded commutative algebra, and $\phi \in
\Fk(\mc M)\otimes R$ satisfies $(\d + \tr)\phi = 0$.  Then $\phi$
carries over to an element of $ \Fk(V_\Tate)\otimes R  $ and so,
when we pass to homology, an element of $\Fk ( \mc H) \otimes R$.
This allows us to include various tautological classes, such as
kappa classes, etc.

As explained in section \ref{section gw invariants}, choice of a
complementary subspace to $\mc H_+$, preserved by $t^{-1}$, leads to
an isomorphism
$$H_\ast(\mf F(\mc H)) \iso \Sym^\ast t^{-1} H_\ast(V)[t^{-1}] $$
and thus a more familiar looking Gromov-Witten potential. Changing
the polarisation changes this potential  by an element of Givental's
twisted loop group.

\begin{remark}
Maxim Kontsevich has informed me that he independently discovered
the relation with Givental's group.
\end{remark} 

\begin{remark}
The main point at which this informal exposition differs from the
rigourous construction contained in the rest of the paper is the
following.  Recall that we don't really have a natural
transformation $C_\ast(\mc M) \to F$, but instead we have a kind of
semi-direct product functor $F^{\mc M}$.  It turns out that we don't
end up with an element in the BV algebra $\Sym^\ast V_{hS^1}$, but
instead we construct a module for the Weyl algebra $\mc
W(V_\Tate)[[\lambda]]$ together with an element in it.  The homology
of this module is a module for $\mc W( \mc H)[[\lambda]]$. Let us
continue to assume that the pairing on $H_\ast(V)$ is
non-degenerate, and that the operator $\op D : V \to V$ is zero on
homology. Then this module is irreducible, and is isomorphic as a
$\mc W(\mc H)[[\lambda]]$ module to $\Fk(\mc H)[[\lambda]]$, in a
unique way up to scale. Thus we find a state in $\Fk(\mc
H)[[\lambda]]$. This is encoded in the ideal in the Weyl algebra
which annihilates it.
\end{remark} 

\subsection{The holomorphic anomaly}

This picture connects very well with Witten's approach to the
holomorphic anomaly \cite{wit1993}.  In Witten's picture, as I
understand it, we interpret the B model potential (without
descendents) for a Calabi-Yau $3$-fold $X$  as an element of the
Fock space associated to the symplectic vector space $H^3(X)$. On
the moduli space of Calabi-Yaus, there is a Gauss-Manin connection
on $H^3(X)$, which preserves the symplectic form.  This therefore
induces a projectively flat connection on the associated Fock space,
and the line spanned by the potential should be flat.

Witten doesn't phrase things in quite this way. Rather, he thinks of
a single Fock space, associated to $H^3(X)$, for some Calabi-Yau
$X$.  This doesn't depend on the complex structure on $X$. Each
choice of complex structure, however, yields a polarisation
$$H^3(X,\C) = ( H^{3,0} \oplus H^{2,1} ) \oplus (H^{1,2} \oplus
H^{0,3})$$ into a direct sum of Lagrangian subspaces.  Also each
complex structure on $X$ yields a B model potential, which can be
considered to be a function on $H^{3,0} \oplus H^{2,1}$.  Using the
polarisation, we can identify the space of functions on $H^{3,0}
\oplus H^{2,1}$ with the Fock space for $H^3(X,\C)$, and the claim
is that the line in the Fock space is independent of the choice of
complex structure on $X$ in a given connected component of the
moduli space of complex structures.

This picture is of course equivalent to the one in the previous
paragraph.

Let us now see how this works in our context. For each Calabi-Yau
$A_\infty$ category, we have a Fock space with an element in it.  We
can think of this as a sheaf of left ideals in the sheaf of Weyl
algebras on the CY moduli space \footnote{ Really, we have an ideal
in something quasi-isomorphic to this Weyl algebra}.

The Weyl algebra is associated to the periodic cyclic chain complex,
shifted by $d$. There should be a Gauss-Manin flat connection on
this, which is the chain level version of the one introduced by
Getzler \cite{get1993} on periodic cyclic homology. The latter was
used in Barannikov's work on the B model \cite{bar1999, bar2000}.

\begin{conjecture}
After modifying the Gromov-Witten potential to take account of the
unstable moduli spaces $(g,n)= (0,1)$ and $(0,2)$, the ideal is
preserved by this flat connection.
\end{conjecture}
This will be discussed elsewhere.

For a Calabi-Yau $X$ over $\K$,  we use an appropriate $A_\infty$
version of the derived category of sheaves.  Then we can use the
Hochschild-Kostant-Rosenberg theorem to identify periodic cyclic
homology with $H^{-\ast}_{DR}(X) ((t))$, where $t$ has degree $-2$.
The Gauss-Manin connection on this is the $\K((t))$ linear extension
of the usual one.

\vspace{7pt} It is interesting to note that the ideal in the Weyl
algebra can be defined even for degenerate TCFTs, where the pairing
on $V$ is degenerate on homology.  The TCFT constructed from a
non-compact symplectic manifold should yield an example of such.
Calabi-Yau $A_\infty$ categories where the pairing on Hochschild
homology is degenerate can be thought of as lying on the boundary of
the moduli space of Calabi-Yau $A_\infty$ categories, corresponding
to large complex structure (B-model) or large volume (A-model)
limits. This idea is made clear in Seidel's work \cite{sei2002},
where the Fukaya category of a projective variety is seen as a
deformation of the Fukaya category of an affine piece.

\subsection{Open closed Gromov-Witten invariants}
In future work I plan to consider the open-closed version of these
constructions. In a similar way to the closed case, Zwiebach
\cite{zwi1998} has constructed a Batalin-Vilkovisky algebra
 from the moduli spaces of Riemann surfaces with open and
closed boundary. Again, there is a unique up to homotopy solution of
the quantum master equation in this BV algebra satisfying  certain
properties, which plays the role of the fundamental chain in these
open-closed moduli spaces.  This corresponds to the fundamental
chain of the moduli space of surfaces with open-closed markings
constructed by Liu \cite{liu2002}, which is an orbifold with
corners.

For a Calabi-Yau $A_\infty$ category,  these fundamental chains give
operators between spaces of morphisms in the category, and the
Hochschild complex.  This structure is a kind of quantisation of the
$A_\infty$ structure.  For the $A$ model, these operators should
correspond to ``counting'' surfaces  with Lagrangian boundary
conditions and with marked points constrained to lie in certain
cycles.

\subsection{Relations with Barannikov's work}
S. Barannikov \cite{bar1999,bar2000} has previously constructed the
genus $0$ B model potential. His construction works for a Calabi-Yau
$A_\infty$ category satisfying the same conditions as are used here.

His idea is that the genus zero potential is encoded in a flat
connection on the periodic cyclic chain complex on the moduli space
of Calabi-Yau $A_\infty$ categories. The periodic cyclic chain
complex is quasi-isomorphic to what we have been calling $V_\Tate$.
For each point of the moduli space we have a subspace $ V^{hS^1}$,
which corresponds to negative cyclic chains. This subspace satisfies
a Griffiths transversality condition, giving what Barannikov calls a
semi-infinite variation of Hodge structure. After choosing a
polarisation, Barannikov shows how to construct a Frobenius manifold
from such a semi-infinite variation of Hodge structure.

A formulation of these ideas closer to what we are doing here has
been given by Givental \cite{giv2001, giv2004}. Fix a reference
Calabi-Yau $A_\infty$ category, $\mc C$. Each nearby category $\mc
C'$ has a subspace $tV^{hS^1}(\mc C') \subset V_\Tate(\mc C')$. We
can translate these to subspaces $tV^{hS^1}(\mc C') \subset
V_\Tate(\mc C)$ using the flat connection on $V_\Tate$. Here they
sweep out a Lagrangian cone.  If we choose a polarisation of
$V_\Tate(\mc C)$, then the cone is the graph of a function on the
positive part, which Givental shows satisfies the equations of a
genus $0$ potential. Givental's twisted loop group acts by change of
polarisation.

Our construction yields an ideal in the Weyl algebra for $V_\Tate$.
The semi-classical limit of this is an ideal in the symmetric
algebra of $V_\Tate$, which cuts out a Lagrangian submanifold in the
dual of $V_\Tate$, which is quasi-isomorphic to a completion of
$V_\Tate$. After taking account of the moduli spaces of curves with
$(g,n) = (0,1)$ and $(0,2)$, this should correspond to the cone
constructed by Barannikov and Givental.

\section{Complexes with a circle action}

Most of the rest of the paper consists of going through this
construction in more detail.   I will give all the definitions of
the previous section again, but more carefully.   Firstly,  we
consider again complexes with a circle action.

Let $V$ be a chain complex, which is either $\Z$ or $\Z/2$ graded. A
circle action on $V$ is by definition an action of the dga
$H_\ast(S^1)$.  This consists of a map $\op{D} : V \to V$, which is
of square zero, commutes with $\op d$, and in the $\Z$ graded case
is of degree $1$.

There are several natural associated complexes, as explained in
\cite{jon1987,hoo_jon1987,lod1998}. The first is the \emph{homotopy
invariants}
$$
V^{hS^1} = V[[t]]
$$
with differential $\op d + \op t D$.   That is if $v f(t)$ is an
element $V[[t]]$, then
$$\op d ( v f(t)) = (\op d v) f(t) + (\op D v) t f(t)$$

The \emph{Tate complex} is
$$
V_{\Tate} = V((t))
$$
again with differential $\op d + t \op D$.

The \emph{homotopy coinvariants} is the space
$$
V_{hS^1} = t^{-1} V[t^{-1}] = V_{\Tate} / V^{hS^1}
$$
with differential induced from that on $V_{\Tate}$.

\begin{remark}
The use of the completed spaces $K[[t]]$ and $K((t))$ is essential.
If we used $K[t]$ and $K[t,t^{-1}]$, then the functors sending $V
\to V^{hS^1}$ or $V_\Tate$ would not be exact.  In the definition of
the various cyclic (co)homology groups it is also essential to
complete in this way.
\end{remark}

If $V$ is graded, and not just $\Z/2$ graded, then all of these
spaces are graded.  The grading on $V^{hS^1}$ and $V_{\Tate}$ is
defined by giving $t$ degree $-2$. The grading on $V_{hS^1}$ is
defined by giving $t^{-k}$ degree $2k-2$.

Let
$$
E_{S^1} = \K((t))[\epsilon]/ \K[[t]][\epsilon] = t^{-1}
\K[t^{-1},\epsilon]
$$
with the differential $\op{d} f = - \epsilon t f$ and the circle
action $\op{D} f = \epsilon f$. Here we give $t^{-k}$ degree $2k-2$
as before, and $\epsilon$ degree one.

Note that the coinvariants for the $H_\ast(S^1)$ action on $V
\otimes E_{S^1}$ is $V_{hS^1}$. That is,
$$
V \otimes_{H_\ast(S^1)} E_{S^1} = V_{hS^1}
$$
Also $E_{S^1}$ is a flat $H_\ast(S^1)$ module, quasi-isomorphic to
$\K$.  So that $V_{hS^1} = V \otimes^{\mbb L}_{H_\ast(S^1)} \K$.

Similar constructions exist when there are $n$ commuting circle
actions, i.e. an action of $H_\ast((S^1)^n)$.  We have to be a
little careful about completions here. The space we use to define
homotopy coinvariants is $\K[[t_1,\ldots,t_n]]$. That used to define
the Tate complex is $K((t_1,\ldots,t_n))$, which by definition
consists of series
$$\sum \lambda_{i_1,\ldots,i_n} t_1^{i_1} \ldots t_n^{i_n}$$
such that  $\lambda_{i_1,\ldots,i_n} = 0$ whenever $\op{min}
(i_1,\ldots,i_n)$ is sufficiently small.

  We are also
interested in complexes with an action of $H_\ast(S^1 \wr S_n)$.  If
$V$ is such a complex, we let
$$
V_{hS^1 \wr S_n} = \left( V  \otimes t_1^{-1} \ldots t_n^{-1}
\K[t_1^{-1}, \ldots, t_n^{-1} ] \right)_{S^n}
$$
where the subscript $S_n$ refers to coinvariants, and the
differential is $\op d + \sum t_i \op D_i$.

\subsection{Relation with the equivariant homology of spaces}
Let $X$ be a (reasonable) topological space with an $S^1$ action.
Let $ES^1$ be a contractible space with an $S^1$ action.
\begin{proposition}
There is a natural isomorphism
$$
H_\ast ( X \times_{S^1} ES^1 ) = H_\ast ( C_\ast(X) _{hS^1})
$$
The action of $\K[[t]]$ on the right hand side corresponds to  cap
product with the first Chern class of the $S^1$ bundle $X \times
ES^1$ over $X \times_{S^1} ES^1$.
\end{proposition}
\begin{proof}
This result is true in much greater generality.  $S^1$ could be
replaced by a topological monoid or topological category satisfying
some mild topological conditions.  The best proof in this generality
involves simplicial methods.  There is a simplicial space model for
the homotopy quotient $X//S^1 \simeq X \times_{S^1} ES^1$. Singular
chains on this give a simplicial object of $\Comp_\K$.  This should
be compared with  a simplicial model for the homotopy tensor product
$C_\ast(X) \otimes_{C_\ast(S^1)}^{\mbb L} \K$.

Instead of going through this argument, I will sketch a geometric
construction of the map $H_\ast(C_\ast(X)_{hS^1}) \to H_\ast(X
\times_{S^1} ES^1)$, which makes it clear that multiplication by $t$
corresponds to cap product with the first Chern class.

Let $S^\infty \subset \C^\infty$ be the set of vectors of norm $1$.
It is well known that $S^\infty$ is contractible, and that the
natural $S^1$ action is free. The quotient $S^\infty/S^1$ is
$\mbb{CP}^\infty$, which is a model for $BS^1$.

Recall that the complex $E_{S^1} = t^{-1} \K[t^{-1},\epsilon]$ has
differential $\op d f = - \epsilon t f$ and circle action $\op D f =
\epsilon f$. The complex $C_\ast(S^\infty)$ also has a circle
action, coming from the map $C_\ast(S^1) \otimes C_\ast(S^\infty)
\to C_\ast(S^\infty)$ and the fundamental chain in $C_\ast(S^1)$. We
define a map $E_{S^1} \to C_\ast(S^\infty)$, compatible with the
differential and the circle action, as follows. We send $t^{-1}$ to
the point $(1,0,\ldots)$, considered as a zero chain. Then $\epsilon
t^{-1}$ must go to the circle $(z,0,\ldots)$ with $\abs {z} = 1$,
with its canonical anti-clockwise orientation. Since $-t^{-2}$
bounds $\epsilon t^{-1}$, we can send $-t^{-2}$ to a fundamental
chain of the cycle $(z_1,z_2,0,0,\ldots..) \in S^\infty$ with $z_2
\in [0,1]$. This is oriented in a canonical way, as it is isomorphic
to the disc $\abs{z_1} \le 1$. We can continue on in this fashion,
and find that $t^{-k}$ gets send to $(-1)^{k+1}$ times a fundamental
chain for $(z_1,z_2,\ldots,z_k,0\ldots)$ with $z_{k} \in [0,1]$, and
$\epsilon t^{-k}$ gets sent to $(-1)^{k+1}$ times a fundamental
chain of the $2k-1$ sphere $(z_1,\ldots,z_k,0,\ldots)$. The sphere
is oriented as the boundary of the ball
$(z_1,\ldots,z_k,0,\ldots)\subset \C^k \subset \C^\infty$ with
$\norm{z} \le 1$.

This map induces a map $C_\ast(X) \otimes_{H_\ast(S^1)} E_{S^1} \to
C_\ast(X \times_{S^1} S^\infty)$, which is a quasi-isomorphism.

Similar remarks hold for cohomology.   To check that cap product by
the first Chern class corresponds to multiplication by $t$, all we
have to check is the sign.  We can do this on $BS^1$. Note that the
line bundle over $\mbb{CP}^\infty$ corresponding to the principal
$S^1$ bundle $S^\infty \to \mbb{CP}^\infty$ is $\Oo(-1)$.  Now,
$t^{-k}$ corresponds to $(-1)^{k-1}[\mbb{CP}^{k-1}] \in
H_{2k-2}(\mbb{CP}^\infty)$, which makes it clear that multiplication
by $t$ is the same as cap product with $c_1(\Oo(-1))$.

\end{proof}

\section{The category of annuli and functors from it}
\label{section functors annuli}

For an integer $m$, let $\mc M_g(m)$ be the moduli space of
connected Riemann surfaces of genus $g$ with $m$ boundary
components, considered to be outgoing. We allow $m = 0$.  Also we
allow ``unstable'' surfaces; the only restrictions are that $g \ge
0$, $m \ge 0$.  However, we need to treat the cases when $g = 0,1$
and $m = 0$ separately. Since we can't glue anything to surfaces in
these spaces, these are essentially placeholders.  We declare that
$\mc M_0(0) = \mc M_1(0)$ are both a point.

Technically, the spaces  $\mc M_g(0)$ are topological stacks.  We
will always work with the coarse moduli space.  This is reasonable
in our setting, as ultimately we only care about rational singular
chains.   We will somewhat loosely use the language of orbispaces.
For instance, we will say that $X$ is a principal $S^1$ orbi-bundle
over $Y$ to mean that $Y$ is the coarse moduli space of an orbispace
over which $X$ is a principal $S^1$ bundle.

The boundary components of the surfaces in $\mc M$ are parameterised
with the opposite orientation to that induced from the orientation
on the surface. That is, if we take the vector field on the boundary
associated to the parameterisation, and apply the complex structure
$J$ to it, it becomes outward-pointing. On $\mc M_g(n)/(S^1)^n$  we
have $n$ principal $S^1$ orbi-bundles. The associated complex line
bundles correspond to the tangent lines at the marked point of a
punctured curve.

Define $\mc M(m)$ like $\mc M_g(m)$ except that the surfaces need
not be connected.  As before, we consider any two complex structures
on a torus or sphere with no boundaries to be the same.

Let
$$
\mc M^s(m) \subset \mc M(m)
$$
be the subspace of stable surfaces, that is those surfaces each of
whose connected components have negative Euler characteristic.

Sometimes it will be more convenient to use $\mc M^s$, and sometimes
$\mc M$.  The main advantage of using $\mc M$ is that it includes
the operation of forgetting a boundary component; if we just used
$\mc M^s$ we would lose information.  On the other hand, using $\mc
M^s$ makes notation much simpler when we compare the solution of the
master equation with the fundamental class of Deligne-Mumford space.

Now we define a topological symmetric monoidal category $A$, which
is a subcategory of $\mc S$.  The objects of $A$ are the
non-negative integers, and the morphisms are the morphisms in $\mc
S(n,m)$ given by Riemann surfaces each of whose connected components
is  an annulus.  As each such annulus has at least one incoming
boundary component, this is a subcategory.  Sending $m \mapsto \mc
M(m)$ defines a symmetric monoidal functor $A \to \op{Top}$.  The
maps $A(m,n) \times \mc M(m) \to \mc M(n)$ are given by gluing
annuli onto the boundary of the surfaces in $\mc M(m)$.

Note that $\mc M^s(m) \subset \mc M(m)$ is a sub-functor.

If we used actual annuli, the categories $\mc S$ and $A$ would not
be unital.  So let us modify the definition a little, to something
weakly equivalent.  In $\mc S$ instead of annuli we now use
``infinitely thin'' annuli, i.e. circles. The parameterisations on
each ``boundary'' of the infinitely thin annulus are then required
to differ from each other only by a rotation and possibly (if both
boundaries are incoming) a change of orientation.   With this
definition, $A$ is a unital category, and $A(1,1) = S^1$. Also
$A(n,n) = S^1 \wr S_n$ as a group.

This modification doesn't change anything essential, as in
\cite{cos_2004oc} I showed that quasi-isomorphic symmetric monoidal
categories have homotopy equivalent categories of functors.

Let $C_* : \op{Top} \to \Comp_\K$ denote the functor of  normalised
singular simplicial chains with $\K$ coefficients. Normalised means
we quotient out by degenerate simplices. This is a symmetric
monoidal functor: the monoidal structure comes from the shuffle
product maps $C_\ast(X) \otimes C_\ast(Y) \to C_\ast(X \times Y)$.

Therefore  we get a functor
\begin{align*}
C_\ast(\mc M) : C_\ast(A) &\to \op{Comp}_{\K} \\
m & \mapsto C_\ast(\mc M(m))
\end{align*}
This is not quite the functor we need, for a technical reason.  Let
$\mc M_{conn}(m) \subset \mc M(m)$ be the subspace of connected
surface. For a finite set $I$, let $\mc M_{conn}(I)$ be the moduli
space of connected surfaces where the boundaries are labelled by the
set $I$. Let $[m] = \{1,\ldots,m\}$. Consider the space
$$
C'_\ast(\mc M(m)) = \oplus_k \left (\oplus_{[m] = I_1 \amalg \ldots
\amalg I_k } C_\ast(\mc M_{conn}(I_1) ) \otimes \ldots C_\ast(\mc
M_{conn} (I_k )) \right) _{S_k}
$$
where the sum ranges over decompositions of the set $[m]$ into a
disjoint union of possibly empty subsets.  It is clear that there is
a natural quasi-isomorphism $C'_\ast(\mc M) \to C_\ast(\mc M)$ and
that $C'_\ast(\mc M)$ defines a functor $A \to \op{Comp}_{\K}$.  We
can think of $C'_\ast(\mc M)$ as the sub-functor of $C_\ast(\mc M)$
generated by connected surfaces.

Let us give a generators and relations description of the category
$H_\ast(A)$, as a unital symmetric monoidal category with objects
$\Z_{\ge 0}$.  Generators are
\begin{enumerate}
\item
The fundamental class $\op D \in H_1(A(1,1))$.
\item
The class of a point in $A(2,0)$, the moduli space of annuli with
two incoming boundaries.  Call this $G \in H_0(A(2,0))$.
\end{enumerate}
Relations are
\begin{align*}
\op D^2 &= 0 \\
G \circ ( \op D \otimes 1 - 1 \otimes \op D ) &= 0
\end{align*}

\begin{lemma}
There is a quasi-isomorphism $H_\ast(A) \to C_\ast(A)$.
\end{lemma}
\begin{proof}
It suffices to write down the map on the generating morphisms of
$H_\ast(A)$.   The morphism $\op D$ goes to a fundamental chain in
$C_1(A(1,1))$.  The morphism $G$ goes to the chain associated to an
annulus in $C_0(A(2,0))$. Pick the annulus where both
parameterisations start at the same point.

It is easy to check the relations hold.  It is crucial here that
we use the normalised singular simplicial chain complex.
\end{proof}

Thus, instead of considering functors from $C_\ast(A)$, we will
always use functors from $H_\ast(A)$.   We can describe such
functors explicitly, using the generators and relations description
for $H_\ast(A)$.  A functor $F : H_\ast(A) \to
\op{Comp}_{\K}^{\Z/2}$ is given by:
\begin{enumerate}
\item
For each $n \ge 0$, a $\Z/2$ graded complex $F(n)$, with maps $F(n)
\otimes F(m) \to F(n+m)$, and  $S_n$ actions on $F(n)$.
\item
For each $1 \le i \le n$, an odd operator $\op D_i : F(n) \to F(n)$.
\item
Maps $G_{ij} : F(n) \to F(n-2)$, for each $1 \le i < j \le n$.
\end{enumerate}
This data satisfies some straightforward axioms, most of which
simply express the fact that the operators $G_{ij}, \op D_i$
interact well with the symmetric group actions and tensor products.
Some other axioms are :
\begin{align*}
\op D_i^2 = [\op D_i, \d] = [\op D_i, \op D_j] &= 0 \\
[G_{ij},\d ] = G_{ij} \circ ( \op D_i - \op D_j) &= 0
\end{align*}

A particularly simple case happens when the functor is \emph{split},
that is the maps $F(1)^{\otimes n} \to F(n)$ are isomorphisms, for
all $n$ (including $n = 0$, when we find $F(0) = \K$). In this case,
let $V = F(1)$.

\begin{lemma}
A split symmetric monoidal functor $F : H_\ast(A) \to
\op{Comp}_\K^{\Z/2}$ is described by:
\begin{enumerate}
\item
A complex $V$.
\item
An odd operator $\op D : V \to V$ which is of square zero and
commutes with $\d$, i.e.\ a circle action.
\item
An even symmetric pairing $\ip{\quad}$ on $V$, such that
\begin{align*}
\ip{\d v,v' } +  (-1)^{\abs{v}}\ip{v, \d v' }  &= 0 \\
\ip{\op D v,v' } -  (-1)^{\abs{v}}\ip{v, \op D v' }  &= 0
\end{align*}
\end{enumerate}
\end{lemma}
In the $\Z$ graded case, the operator $\op D$ is of degree one.

\subsection{The functors associated to a TCFT}
A TCFT is a functor $F : C_\ast(\mc S) \to \op{Comp}_{\K}^{\Z/2}$,
which is h-split, i.e. the maps $F(1)^{\otimes n} \to F(n)$ are
quasi-isomorphisms. In particular, $F$ restricts to a functor
$H_\ast(A) \to \op{Comp}_{\K}^{\Z/2}$. This functor $H_\ast(A) \to
\Comp_\K^{\Z/2}$ associated to a TCFT only encodes a very small
amount of the structure of the TCFT.  One could hope that there
would be a natural transformation $C'_\ast(\mc M) \to F$ of functors
on $H_\ast(A)$. However, because of the restriction that the
morphism surfaces in $\mc S$ have at least one incoming boundary,
this is not in general true.

Instead, we will construct a semi-direct product functor $F^{\mc M}$
which encodes the data of the action of $C_\ast(\mc S)$ on $F$.

For an integer $m$, we define
$$
F^{\mc M} (m) = \oplus_{[m] = I \amalg J} F(I) \otimes C'_\ast(\mc
M(J))
$$
where the direct sum is over decompositions of the set $[m]$ into
two possibly empty disjoint subsets.

To define the structure of functor from $H_\ast(A)$, we need to
write down gluing maps $G_{ij} : F^{\mc M}(m) \to F^{\mc M}(m-2)$,
and circle actions $\op D_i : F^{\mc M}(m) \to F^{\mc M}(m)$.

On the direct summand $F(I) \otimes C'_\ast(\mc M(J))$, the circle
action $\op D_i$ is the corresponding one on $F(i)$ or $C'_\ast(\mc
M(J))$, depending on whether $i \in I$ or $i \in J$.

We will define the gluing maps also on the direct summand $F(I)
\otimes C'_\ast(\mc M(J))$. If $i,j \in I$, or $i,j \in J$, then
$G_{ij}$ is the gluing map from $F$ or $C'_\ast(\mc M)$.  If $i \in
I$ and $j \in J$, then the gluing map is more difficult to
construct.  This uses the action of $C_\ast(\mc S)$ on $F$. It is
enough to define this map on connected surfaces, $C_\ast(\mc
M_{conn}(J))$.  Then we can turn the $j$ boundary around to give an
element in $C_\ast(\mc S(j,J\setminus \{j\})$ which acts on $F$.

More formally, by the definition of $C'_\ast(\mc M)$, we have
$$
F^{\mc M} (m) = \oplus_k \left( \oplus_{[m] = I \amalg J_1 \amalg
\ldots J_k} F(I) \otimes C_\ast(\mc M_{conn}(J_1)) \otimes
 \ldots C_\ast(\mc M_{conn}(J_k))
\right)_{S_k}
$$
Now let $i \in I$ and $j \in J_1$; we will define the gluing map on
one of the direct factors of this decomposition.  For simplicity, we
will assume $k = 1$.

Then there is an isomorphism
$$
C_\ast(\mc M_{conn}(J)) \iso C_\ast(\mc S_{conn}(j, J\setminus\{j\})
= C_\ast(\mc S_{conn}(i, J\setminus\{j\}))
$$
There is a map
$$ C_\ast(\mc S_{conn}(i, J\setminus\{j\})) \otimes C_\ast( \mc S (I \setminus\{i\}, I \setminus\{i\} ) )
\to  C_\ast(\mc S( I , J\setminus\{j\} \amalg I \setminus\{i\} ))$$
coming from the symmetric monoidal structure on $\mc S$, given by
disjoint union.

Placing the identity morphism $I \setminus\{i\} \to I
\setminus\{i\}$ on the second factor gives a map
$$
C_\ast(\mc S_{conn}(i, J\setminus\{j\})   \to  C_\ast(\mc S( I ,
J\setminus\{j\} \amalg I \setminus\{i\} )
$$
$F$ is a functor $C_\ast(\mc S) \to \op{Comp}_{\K}^{\Z/2}$. There is
an action map
$$
 C_\ast(\mc S( I ,
J\setminus\{j\} \amalg I \setminus\{i\} ) ) \otimes F(I) \to F(
J\setminus\{j\} \amalg I \setminus\{i\} )
$$
Composing these maps gives the required operator
$$
G_{ij} : F(I) \otimes C_\ast(\mc M_{conn}(J)) = F(I) \otimes
C_\ast(\mc S_{conn}(i, J\setminus\{j\})) \to F( J\setminus\{j\}
\amalg I \setminus\{i\} )
$$
It is not difficult to check that this defines a functor $H_\ast(A)
\to \op{Comp}_\K^{\Z/2}$.

\section{The Weyl algebra and the Fock space associated to a
functor}

 \label{section weyl fock} Let $F : H_\ast(A) \to
\op{Comp}_\K^{\Z/2}$ be a symmetric monoidal functor.  We will
construct an associated Weyl algebra and Fock space.

\subsection{ The construction in a simplified case.}

 Let us first consider
the simplified case when $F$ is split. Let $V = F(1)$. Then $V$ has
a circle action, and  we have the auxiliary equivariant chain
complexes, $V_{hS^1}$, $V^{hS^1}$ and $V_{\Tate}$.

Define an antisymmetric form $\Omega$ on $V_{\Tate}$ by
$$
\Omega ( v f(t), w g(t) ) = \ip{v,w} \op {Res} f(-t) g(t) \op d t
$$
This is the same as the form used in the work of Givental and Coates
\cite{coa_giv2001,giv2001, giv2004}.  In the case when the inner
product on $V$ is non-degenerate this is symplectic.  Note that
$\Omega$ is compatible with the differential, that is
$$
\Omega ( \op d (v f(t)), w g(t) )  + (-1)^{\abs v} \Omega ( v
f(t), \op d ( w g(t)) ) =0
$$
This follows from the fact that on $V$, $\op d$ is skew self adjoint
and $\op D$ is self adjoint with respect to the pairing
$\ip{\quad}$. Thus, we have an associated Weyl algebra $\mc
W(V_\Tate)$. $V_\Tate$ is polarised, as $V_\Tate = V_{hS^1} \oplus
V^{hS^1}$. The differential on $V_\Tate$ preserves $V^{hS^1}$, but
not in general $V_{hS^1}$. Let $\Fk(V_\Tate)$ be the associated Fock
space. This is defined to be the quotient of $\mc W(V_\Tate)$ by the
left ideal generated by $V^{hS^1}$.    We can identify
$$
\Fk(V_\Tate) = \Sym^\ast V_{hS^1}
$$
As, we can consider $\Sym^\ast V_{hS^1}$ as a subalgebra of $\mc
W(V)$, using the splitting of the map $V_{\Tate} \to V_{hS^1}$. Then
the action on this on the element $1 \in \Fk(V_\Tate)$ gives the
isomorphism.  This is not an isomorphism of complexes, however,
because $V_{hS^1} \subset V_\Tate$ is not a subcomplex.

Let us write the natural differential on $\Fk(V_\Tate)$ as
$\what{\op d}$. This is the differential obtained by realising it as
a quotient of $\mc W(V_\Tate)$. This is an order $2$ differential
operator.   Let $\op d$ denote the usual differential on $\Sym^\ast
V_{hS^1} $, which we identify with $\Fk(V_\Tate)$. Then we can write
$$
\what{\op d} = \op d +  \tr
$$
where $\tr$ is an odd order 2 differential operator on
$\Fk(V_\Tate)$,   and satisfies
$$
[\op d,\tr] = \tr^2 = 0
$$
We can describe $\tr$ explicitly.  It is an order $2$ differential
operator on $\Sym^\ast V_{hS^1}$. Such an operator is uniquely
characterised by its behaviour on $\Sym^{\le 2} V_{hS^1}$.  $\tr$ is
zero on $\Sym^{\le 1} V_{hS^1}$, and for $(v_1 f_1(t_1))(v_2
f_2(t_2)) \in \Sym^2 V_{hS^1}$, we have
$$
\tr((v_1 f_1(t_1))(v_2 f_2(t_2)) ) = \ip{\op D v_1, v_2 } \op{Res}
f_1(t_1) f_2(t_2) \op d t_1 \op d t_2
$$
This has been proved in lemma \ref{lemma quantised bv}.
\subsection{The construction in general}
We want to mimic this construction in general. Let $F : H_\ast(A)
\to \op{Comp}_\K^{\Z/2}$ be any symmetric monoidal functor.  On
$F(n)$ there are $n$ commuting circle actions, that is operators
$\op D_i$ for $1 \le i \le n$, which (super)-commute and square to
zero.  Thus we can form the various auxiliary complexes,
\begin{align*}
F_{\Tate}(n) &= F(n) \otimes \K((t_1,\ldots,t_n)) \\
F^{hS^1}(n) &= F(n) \otimes \K[[t_1,\ldots,t_n]] \\
F_{hS^1}(n) &= F(n) \otimes \K((t_1))/\K[[t_1]] \otimes \ldots
\K((t_n))/\K[[t_n]]
\end{align*}
with differential
$$
\op d + \sum t_i \op D_i
$$
Let $G_{ij} : F(n) \to F(n-2)$ be the gluing map, coming from the
class of a point in $H_0(A(2,0))$.  For $1 \le i < n$ denote by
$$\Omega_{i} : F_\Tate(n) \to F_\Tate(n-2)$$
the map defined by
$$
a \otimes f(t_1,\ldots,t_n) \mapsto G_{i,i+1}(a) \otimes \op{Res}_{z
=  0} f(t_1,\ldots, t_{i-1},  -z, z, t_i,\ldots, t_{n-2}) \op dz
$$
For each $1 \le i <  n$, let $\sigma_{i} \in S_n$  be the
transposition of $i$ with $i+1$.  Recall $S_n$ acts on $F(n)$; this
action extends to each of the auxiliary complexes mentioned above.

There are tensor product maps $F_\Tate(n) \otimes F_\Tate(m) \to
F_\Tate(n+m)$, and similarly for $F_{hS^1}$ and $F^{hS^1}$.  The
space $\oplus_n F_\Tate(n)$ is an associative algebra, with product
coming from these tensor product maps.

We define the \emph{Weyl algebra} $\mc W(F)$ to be the quotient of
$\oplus_n F_\Tate(n)$  by the two-sided ideal generated by the
relation,
$$
x - \sigma_i(x) =  \Omega_i(x)
$$
for each $x \in F_\Tate(n)$.

The \emph{Fock space} $\mf F(F)$ is defined to be the quotient of
$\mc W(F)$ by the left ideal spanned by those elements $x \in
F_{\Tate}(n)$ which contain no negative powers of $t_n$.

We can consider
$$
\oplus F_{hS^1}(n)_{S_n}
$$
to be a subalgebra of $\mc W(F)$, using the standard splitting of
the map $F_\Tate(n) \to F_{hS^1}(n)$.  Here the subscript $S_n$
refers to coinvariants, so that $\oplus F_{hS^1}(n)_{S_n}$ is a
commutative algebra. The action of $\oplus F_{hS^1}(n)_{S_n}$ on the
vector $1 \in \mf F(F)$ generates $\mf F(F)$, and induces an
isomorphism
$$
\mf F (F) \iso \oplus F_{hS^1}(n)_{S_n}
$$
As before, this is \emph{not} an isomorphism of complexes.  We will
refer to the natural differential on the left hand side as
$\what{\op d}$, and that on the right hand side as $\op d$. The
differential $\what{\op d}$  is an order $2$ differential operator,
whereas $\op d$ is a derivation.

It is easy to see that, as before,
$$
\tr \defeq \what{\op d} - \op d
$$
is an order two operator which satisfies $\tr^2 = [\d,\tr]= 0$.  As
before, we can write this operator down explicitly.  For $1 \le i <
j \le n$, define a map $\tr_{ij} : F_{hS^1}(n) \to F_{hS^1}(n-2)$ by

\begin{multline*}
\tr_{ij}(a \otimes f(t_1,\ldots,t_n)) = G_{ij}(\op D_i a) \otimes
\\\op{Res}_{z = 0} \op{Res}_{w = 0}
f(t_1,\ldots,t_{i-1},z,t_{i},\ldots,t_{j-2},w,t_{j-1},\ldots,
t_{n-2}) \op d z \op d w
\end{multline*}
(In this expression $z$ is in the $i$'th position and $w$ is in
the $j$'th position, and the remaining places are filled with
$t_1,\ldots,t_{n-2}$ in increasing order).

Then, $\sum_{i < j} \tr_{ij}$ commutes with symmetric group
actions, and so descends to give a map
$$\tr :  \oplus F_{hS^1}(n)_{S_n} \to
\oplus F_{hS^1}(n)_{S_n}$$ It is now not difficult to check that
$\what{\op d } = \op d +  \tr$.  The proof is the same as that of
lemma \ref{lemma quantised bv}, which proves this in the special
case that $F$ is split.

We will need these constructions when $F$ is the functor
$C'_\ast(\mc M)$ associated to moduli spaces of curves.  In that
case we use the notation $\mc W(\mc M)$, $\mf F(\mc M)$.

Note that if $F \to G$ is a natural transformation of functors
$H_\ast(A) \to \op{Comp}_{\K}^{\Z/2}$, there is an associated
homomorphism of Weyl algebras $\mc W(F) \to \mc W(G)$, and a map
$\Fk(F) \to \Fk(G)$ of $\mc W(F)$ modules.

\subsection{Geometric interpretation of the differential on $\Fk(\mc
M)$ }

We know that
$$H_\ast(  C_\ast (\mc M(m) )_{h(S^1)^m}) = H_\ast(\mc
M(m) / (S^1)^m ).$$  In calculating the differential on $\Fk(\mc M)$
we used operators
$$\tr_{ij} : H_\ast(\mc M(m+2)/(S^1)^{m+2})\to H_\ast(\mc M(m) / (S^1)^m)$$
These operators have a geometric interpretation, which gives a
geometric interpretation to the order two part $\tr$ of the
differential $\dhat = \d + \tr$.

Let
$$
Y = \mc M(m+2) / (S^1)^m \times S^1
$$
where we quotient by the circle actions on all boundaries except the
$i$ and $j$ ones, and by the anti-diagonal circle action from the
$i,j$ boundaries.

There is a map $\pi : Y \to \mc M(m+2)/(S^1)^{m+2}$, whose fibres
are oriented circle bundles.  Also there is a gluing map $\iota : Y
\to \mc M(m)/(S^1)^m$.
\begin{lemma}
The map $\tr_{ij} : H_\ast(\mc M(m+2)/(S^1)^{m+2}) \to H_{\ast + 1}
( \mc M(m)/(S^1)^m)$ is $\iota_\ast \pi^\ast$.
\end{lemma}
\begin{proof}
This consists of unravelling the definition. For simplicity we will
consider the case $m = 0$.

Recall $E_{S^1} = t^{-1} \K[t^{-1},\epsilon]$ is a certain
contractible complex with a circle action.  Then $C_\ast(\mc
M(2))_{h(S^1)^2}$ is the $H_\ast((S^1)^2)$ coinvariants of
$C_\ast(\mc M(2)) \otimes (E_{S^1})^{\otimes 2}$. The latter is
quasi-isomorphic to $C_\ast(\mc M(2))$. The quasi-isomorphism is the
map $C_\ast(\mc M(2)) \otimes (E_{S^1})^{\otimes 2} \to C_\ast(\mc
M(2))$  which sends
$$v \otimes f(t_1,\epsilon_1,t_2,\epsilon_2) \to v \op{Res}
f(t_1,0,t_2,0)$$ This is a map of dg $H_\ast(S^1 \wr S_2)$ modules.

The gluing map $G  : C_\ast(\mc M(2)) \to C_\ast(\mc M(0))$ then
extends to $C_\ast(\mc M(2)) \otimes (E_{S^1})^{\otimes 2}$, by
composition with this quasi-isomorphism.

The map
\begin{align*}
\tr_{12}  : C_\ast(\mc M(2))_{h(S^1)^2} &\to C_\ast(\mc M(0)) \\
v \otimes f (t_1,t_2) & \mapsto G( \op D_1 v ) \op{Res} f(t_1,t_2)
\d t_1 \d t_2
\end{align*}
has an interpretation as : take an element of $C_\ast(\mc
M(2))_{h(S^1)^2}$, lift (in any way) to $C_\ast(\mc M(2)) \otimes
(E_{S^1})^{\otimes 2}$, apply $\op D_1$, then apply the gluing map
there.

This makes the result clear.

\end{proof}

\section{Batalin-Vilkovisky algebras}
\label{section bv}
\begin{definition}
A Batalin-Vilkovisky (BV) algebra is a differential $\Z/2$ graded
super-commutative algebra $B$, together with an odd operator $\tr :
B \to B$, which is an order $2$ differential operator, and satisfies
$$
\tr^2 = [\d, \tr] = \tr(1) = 0
$$
We let $\dhat = \d + \tr$.
\end{definition}
For each functor $F : H_\ast(A) \to \op{Comp}_\K^{\Z/2}$, the Fock
space $\Fk(F)$ constructed in the previous section is a
Batalin-Vilkovisky algebra.

If $B$ is a BV algebra, then it acquires an odd Poisson structure.
The bracket is defined by
$$
\{f,g\} =  \dhat (fg) - (-1)^{\abs {f}} f \dhat(g) - \dhat(f) g
$$
This satisfies the Jacobi identity; see \cite{get1994}.  Here it is
also shown that $\dhat$ is a derivation of this bracket, that is
$$
\dhat \{f,g\} = \{\dhat f,g\} + (-1)^{\abs {f}}\{f,\dhat g\}
$$
Therefore $B$ becomes a differential $\Z/2$ graded Lie algebra, with
this Lie bracket and differential $\dhat$.

The Maurer-Cartan equation in $B$ is the equation
$$
\dhat S + \frac{1}{2} \{S,S\} = 0
$$
This is equivalent to the quantum BV master equation
$$
\dhat \exp(S ) = 0
$$
(whenever this expression makes sense in the algebra $B$). Indeed,
it is easy to see that
$$
\exp(-S) \dhat \exp(S) = \dhat S + \frac{1}{2} \{S,S\}
$$

\subsection{Homotopies between solutions of the master equation}
Consider the differential graded algebra $\K[t,\epsilon]$, where $t$
is of degree $0$ and $\epsilon$ is of degree $-1$, with differential
$\epsilon \frac{\op{d}}{\op{d}t}$.    Let $\mf g$ be a differential
graded Lie algebra, with differential of degree $-1$.  A solution of
the Maurer-Cartan equation in $\mf g$ is an element $S \in \mf
g_{-1}$ satisfying
$$
\op{d} S + \frac{1}{2}[S,S] = 0
$$
If $\mf g$ is only $\Z/2$ graded, $S$ must simply be odd.

A homotopy between solutions $S_0,S_1$ of the Maurer-Cartan equation
in $\mf g$ is an element
$$
S(t,\epsilon) \in \mf g [t,\epsilon]
$$
which satisfies the Maurer-Cartan equation:
$$
\op{d} S + \epsilon \frac{\op{d} S}{\op{d}t} + \frac{1}{2}[S,S] = 0
$$
and such that $S(0,0) = S_0$, and $S(1,0) = S_1$.

Note that we can write
$$
S(t,\epsilon) = S_a(t) + \epsilon S_b(t)
$$
The Maurer-Cartan equation for $S$ implies that $S_a$ satisfies the
Maurer-Cartan equation, and that
$$
\frac{\op{d}S_a(t)}{\op{d}t} = - [ S_b(t), S_a(t)]  - \d S_b(t)
$$
so that the path in $\mf g_{-1}$ given by $S_a(t)$ is tangent to the
action of $\mf g_0$ on solutions of Maurer-Cartan in $\mf g_{-1}$.

Let $\op{MC}(\mf g)$ be the set of Maurer-Cartan elements in $\mf g$
and let $\pi_0(\op{MC}(\mf g))$  be the quotient of this by the
equivalence relation generated by homotopy.  These definitions work
in the $\Z/2$ graded case, and also for odd Lie algebras, with
obvious changes.

If $B$ is a BV algebra, let $\op{BV}(B)$ be the set of solutions of
the master equation in $B$, that is the set of solutions of the
Maurer-Cartan equation in $B$ considered as an odd dg Lie algebra.
Let $\pi_0 \op{BV}(B)$ be the set of homotopy classes of solutions
of the master equation, defined as above.

There is an obvious notion of homotopy between maps $f_0,f_1 : \mf g
\to \mf g'$ of dg Lie algebras. This is a map $F : \mf g \to \mf
g'[t,\epsilon]$ of dg Lie algebras, such that $F(0,0) = f_0$ and
$F(1,0) = f_1$.  Clearly homotopic maps induce the same map  $\pi_0
\op{MC}(\mf g) \to \pi_0\op{MC}(\mf g')$.  It follows that a
homotopy equivalence $\mf g\to \mf g'$ (i.e. a map which has an
inverse up to homotopy) induces an isomorphism on $\pi_0 \op{MC}$.

In nice cases, quasi-isomorphisms of dg Lie algebras also induce
isomorphisms on the set of homotopy classes of solutions of the
Maurer-Cartan equation.  Suppose $\mf g$ is a dg Lie algebra with a
filtration $\mf g = F^1 \mf g \supset F^2 \mf g\supset\ldots$, such
that $\mf g$ is complete with respect to the filtration, and such
that $[F^i \mf g, F^j \mf g] \subset F^{i+j} \mf g$.  In particular
$\mf g /F^2 \mf g$ is Abelian and each $\mf g / F^i \mf g$ is
nilpotent. Then we say $\mf g$ is a filtered pro-nilpotent Lie
algebra.
\begin{lemma}
Let $\mf g,\mf g'$ be filtered pro-nilpotent dg Lie algebras, and
let $f : \mf g \to \mf g'$ be a filtration preserving map. Suppose
the map $\op{Gr} \mf g \to \op{Gr} \mf g'$ induces an isomorphism on
$H_i$ for $i = 0,-1,-2$. Then the map
$$
\pi_0 \op{MC}(\mf g) \to \pi_0 \op{MC} (\mf g')
$$
is an isomorphism. \label{lemma mc isomorphism}
\end{lemma}
\begin{proof}
This result seems to be well known.  For instance it is essentially
theorem 5.1 of \cite{ham_laz2004}, or theorem 2.1 of \cite{get2002}.

For completeness I will sketch a proof.  First we will show that we
can replace $f$ by a surjective map. Let $\mf g'' \subset \mf g
\oplus  \mf g'[t,\epsilon]$ be the subset of elements
$(\gamma,\alpha(t,\epsilon))$ such that $f(\gamma) = \alpha(0,0)$.
This is an analog of the Serre construction which replaces any map
of topological spaces by a fibration.  It is easy to see, by
mimicking the corresponding topological argument, that the natural
maps $\mf g \into \mf g''$ and $\mf g'' \to \mf g$ are inverse
homotopy equivalence.  The key point in the topological argument is
to use the multiplication map $[0,1] \times [0,1] \to [0,1]$.  Here
we instead use a bi-algebra structure on $\K[t,\epsilon]$. The
coproduct is defined on the generators by $t \mapsto t \otimes t$,
$\epsilon \mapsto t \otimes \epsilon + \epsilon \otimes t$.  It is
easy to check that this coproduct is compatible with the
differential.

It remains to show that the map $\pi_0\op{MC}(\mf g'') \to \pi_0
\op{MC}(\mf g')$ is an isomorphism.  There are three things to
prove.
\begin{enumerate}
\item
The map $\op{MC}(\mf g'') \to \op{MC}(\mf g')$ is surjective.
\item
Any two points $T_1,T_2 \in \op{MC}(\mf g'')$ with the same image in
$\op{MC}(\mf g')$ are homotopic.
\item
The map $\op{MC}(\mf g'') \to \op{MC}(\mf g')$ has the path lifting
property.
\end{enumerate}
All of these are proved by the same kind of inductive argument.
\end{proof}
A similar result holds in the $\Z/2$ graded case, under the
assumption that $\op{Gr} f$ is a quasi-isomorphism.

\section{The master equation and the fundamental chain}

\label{section fundamental master}

We have constructed the Sen-Zwiebach Batalin-Vilkovisky algebra
$\Fk(\mc M)$ associated to moduli spaces.  Now we will construct in
this a solution $S$ of the master equation, which  plays the role of
the fundamental class.

There is a natural inclusion map
$$
C_\ast(\mc M_{g}(n))_{hS^1 \wr S_n} \to \Fk(\mc M)
$$
where $\mc M_g(n) \subset \mc M(n)$ is the subspace of connected
surfaces of genus $g$.  Denote by $\Fk^{g,n}(\mc M)$ this subspace.

$\Fk(\mc M)$ is freely generated as a commutative  algebra by the
subspaces $\Fk^{g,n}(\mc M)$.

\begin{proposition}
For each $g,n$ with $2g-2+n > 0$, there exists an element $S_{g,n}
\in \Fk^{g,n} (\mc M )$ of degree $6g-6+2n$, with the following
properties.
\begin{enumerate}
\item
 $S_{0,3}$ is a  $0$-chain of degree $1/3!$ in the moduli
space of Riemann surfaces with $3$ unparameterised, unordered
boundaries.
\item
Form the generating function
$$
S = \sum_{\substack{ g,n \ge 0 \\ 2g-2+n > 0}} S_{g,n}
\lambda^{2g-2+n}
 \in \lambda \Fk(\mc M)[[\lambda]]
$$

$S$ satisfies the Batalin-Vilkovisky quantum master equation :
$$
\dhat e^{S} = 0
$$
Equivalently,
$$
 \dhat{S} + \frac{1 }{2}\{S,S\}  = 0
$$
\end{enumerate}
Further, such an $S$ is unique up to homotopy through such elements.
\label{prop master moduli}
\end{proposition}

A homotopy of such elements is a solution of the master equation in
$\Fk(\mc M) \otimes \K[t,\epsilon]$, satisfying the analogous
conditions.  Here $t$ has degree $0$ and $\epsilon$ has degree $-1$,
and $\op d t = \epsilon$.

\begin{proof}
Let $\mc M_{g,n}$ be the usual moduli space of smooth algebraic
curves of genus $g$ with $n$ marked points. This is rationally
homotopy equivalent to $\mc M_g(n)/(S^1)^n$.

We will need the following bound on the homological dimension of
$\mc M_{g,n}/S_n$:
\begin{equation}
H_i ( \mc M_{g,n} /S_n) = 0 \text{ for } i \ge  6g - 7 + 2n \text{
if } (g,n) \neq (0,3) \label{equation hom_dimn}
\end{equation}
To see this, observe that $\cmod_{g,n}$ is simply connected as an
orbifold, because the mapping class group is generated by Dehn
twists, and compactifying $\mc M_{g,n}$ has the effect of
trivialising the elements of $\pi_1(\mc M_{g,n})$ coming from Dehn
twists. In particular $H_1(\cmod_{g,n}) = 0$. It follows that
$H_1(\cmod_{g,n}/S_n) = 0$, as we are using coefficients in $\K
\supset \Q$.  The boundary of $\cmod_{g,n}/S_n$ is always connected.
(When $(g,n) \neq (0,4)$, the boundary of $\cmod_{g,n}$ is
connected). Poincar\'e duality and the cohomology long exact
sequence for the pair $(\cmod_{g,n}/S_n,\partial \cmod_{g,n}/S_n)$
gives the required bound.

Alternatively, we could use the bounds on the homological dimension
of $\mc M_{g,n}$ obtained by Harer \cite{har1986}.  This gives the
result when $(g,n) \neq (0,4)$.  In that case, it is easy to see
that the coinvariants of the $S_4$ action on $H_1(\mc M_{0,4}) =
\K^2$ are trivial.

Now define a dg Lie algebra $\mf{g}$.  The space $\mf{g}_i$ is the
set of
$$
S = \sum S_{g,n} \lambda^{2g-2+n}  \in \lambda \Fk(\mc M)[[\lambda]]
$$
such that $S_{g,n} \in \Fk^{g,n}(\mc M)$, and $S_{g,n}$ is of degree
$6g-6+2n + 1 + i$.  The bracket $[\quad]_{\mf g}$ is $\{\quad\}$ and
the differential is $\op{d}_{\mf g} =   \dhat$.

The set of homotopy equivalence classes of solutions of the
Maurer-Cartan equation in $\mf{g}$ is the same as the set of
homotopy equivalence classes of solutions $S$ of the master equation
in $\Fk(\mc M)$ satisfying $S_{g,n} \in
 \Fk^{g,n}(\mc M)$ and $S_{g,n}$ is of degree $6g-6+2n$. Filter $\mf g$ by saying $F^k \mf g$
is the set of those $S$ such that $S_{g,n}$ is zero for $2g-2+n <
k$. Then
$$\mf g = F^1 \mf g \supset F^2 \mf g \ldots$$
is a descending filtration by dgla ideals. $\mf g$ is complete with
respect to this filtration.

The bounds (\ref{equation hom_dimn}) on the homological dimensions
of moduli spaces, together with the fact that $\mc M_{0,3}$ is a
point, tell us that
\begin{align*}
H_i ( F^k \mf g / F^{k+1} \mf g ) &= 0  \text{ for } i \ge 0, i = -2 \\
H_{-1} ( F^k \mf g / F^{k+1} \mf g) &=  0 \text { for } k \ge 2 \\
 H_{-1} (\mf g / F^2 \mf g) &= \K
\end{align*}
Therefore the map $\mf g \to \mf g / F^2 \mf g$ satisfies the
conditions of lemma \ref{lemma mc isomorphism}. The result follows
immediately.
\end{proof}

Note that $S$ comes from the stable moduli spaces.  Recall that $\mc
M^s(m)$, the space of surfaces each of whose components has negative
Euler characteristic, is a sub-functor of $\mc M(m)$.  Therefore we
have an associated BV algebra $\Fk(\mc M^s)$ and an injective map of
BV algebras $\Fk(\mc M^s) \to \Fk(\mc M)$. $S$ is in $\Fk(\mc M^s)$.

We want to compare this solution of the master equation with the
usual fundamental class of Deligne-Mumford space. Let $\cmod_{g,n}$
be the space of stable nodal curves of genus $g$ with $n$ marked
smooth points.    Let $\cmod(n)$ be the moduli space of possibly
disconnected stable nodal curves with $n$ marked points.

Define
$$
\Fk(\cmod) = \oplus_n  C_\ast(\cmod(n)/S_n)
$$
This forms a Batalin-Vilkovisky algebra, with BV operator $\tr = 0$,
so that $\dhat = \d$. Let
$$
[\cmod]= \sum \lambda^{2g-2+n} [\cmod_{g,n}/S_n] \in
\Fk(\cmod)[[\lambda]]
$$
where $[\cmod_{g,n}/S_n]$ is a fundamental chain for
$\cmod_{g,n}/S_n \subset \cmod(n)/S_n$.

Let $B$ be a BV algebra.  Then $\lambda B[[\lambda]]$ is a
pro-nilpotent odd Lie algebra, with Lie bracket $\{ \ ,  \ \}$ and
differential $\dhat$.  Therefore we have a set $\op{MC} (\lambda
B[[\lambda]])$ of solutions of the Maurer-Cartan equation in
$\lambda B[[\lambda]]$, or equivalently solutions of the quantum
master equation; and a set $\pi_0 (\op{MC} (\lambda B[[\lambda]]))$
of homotopy classes of solutions. In particular, $S \in \lambda
\Fk(\mc M^s) [[\lambda]]$ and $[\cmod] \in \lambda
\Fk(\cmod)[[\lambda]]$ are such solutions of the Maurer-Cartan
equation.

If $B \to B'$ induces an isomorphism on $H_\ast(B, \dhat) \to
H_\ast(B',\dhat)$, then it induces an isomorphism
$$
\pi_0(\op{MC} ( \lambda B[[\lambda]] )) \simeq \pi_0 ( \op{MC}
(\lambda B'[[\lambda]]))
$$
This follows from lemma \ref{lemma mc isomorphism}.

We say a map in the homotopy category of BV algebras $B \to B'$ is a
map $B'' \to B'$, where $B$ and $B''$ are connected by a sequence of
maps of BV algebras which induce an isomorphism on $\dhat$ homology.
Any such map induces a map
$$
\pi_0(\op{MC} ( \lambda B[[\lambda]] )) \to \pi_0 ( \op{MC} (\lambda
B'[[\lambda]]))
$$
\begin{theorem}
There is a map $\Fk(\mc M^s) \to \Fk(\cmod)$ in the homotopy
category of BV algebras which maps the class $S \in \pi_0(\op{MC} (
\lambda \Fk(\mc M^S)[[\lambda]] ))$ to $[\cmod] \in \pi_0(\op{MC} (
\lambda \Fk(\cmod)[[\lambda]] )) $.
\end{theorem}
\begin{remark}
In fact, we could use $\mc M$ instead $\mc M^s$, but this would
involve some messing around with unstable surfaces.
\end{remark}

We will construct this using a nice model for the spaces $\mc
M^s(m)$, introduced by Kimura, Stasheff and Voronov
\cite{kim_sta_vor1995}.  Let $\mc N_g(n)$ be the moduli space of
algebraic curves in $\cmod_{g,n}$, together with at each marked
point, a ray in the tangent space, and at each node, a ray in the
tensor product of the tangent spaces at each side.   $\mc N_g(n)$ is
a torus bundle over a certain real blow-up of $\cmod_{g,n}$, and is
an orbifold with corners, whose boundary consists of the locus of
singular curves.

Let $\mc N(n)$ be defined in the same fashion, except using possibly
disconnected curves.

The space $\mc N(n)$ has an action of $S^1 \wr S_n$. Also there are
gluing maps $G_{ij} : \mc N(n) \to \mc N(n-2)$, for each $1 \le i <
j \le n$. These satisfy various compatibility conditions, which
means that sending $n \mapsto \mc N(n)$ defines a symmetric monoidal
functor
$$
\mc N : A \to \op{Top}
$$

Passing to chain level, we find a functor $C_\ast(\mc N) : H_\ast(A)
\to \op{Comp}_{\K}$.  Therefore, we have an associated BV algebra
$\Fk(\mc N)$.  Proposition \ref{prop master moduli} applies without
change to $\Fk(\mc N)$, giving a solution of the master equation,
also denoted by $S$, in $\Fk(\mc N)$.

We want to compare this with the solution in $\Fk(\mc M^s)$.  First
we need to compare the two BV algebras. The BV algebra is associated
functorially to a functor $H_\ast(A) \to \op{Comp}_{\K}$.

A natural transformation $F \to G$ between such functors is a
quasi-isomorphism if it induces a quasi-isomorphism on the chain
complexes $F(n) \to G(n)$ for all $n \in \Z_{\ge 0} = \op{Ob}
H_\ast(A)$.  Two functors are quasi-isomorphic if they can be
connected by a chain of quasi-isomorphisms.

If $F,G$ are quasi-isomorphic functors, then the associated BV
algebras $\mf F(F), \mf F(G)$ are quasi-isomorphic on $\op d $
homology (but not necessarily on $\what {\op{d}}$ homology).

\begin{lemma}
The functors $C_\ast(\mc N)$, $C_\ast(\mc M^s)$ are
quasi-isomorphic.
\end{lemma}
\begin{proof}[Sketch of proof]
We will show the corresponding result at the level of the functors
$\mc N, \mc M^s : A \to \op{Top}$.  That is, we will show these
functors are rationally weakly equivalent.  A rational weak
equivalence is a natural transformation that induces an isomorphism
on the rational homology of the associated spaces.

We need to construct  a chain of rational weak equivalences between
$\mc M^s$ and $\mc N$ in the category of functors $A \to \op{Top}$.

Let me sketch such a construction. Firstly, consider a moduli
space $\mc P$ like $\mc N$, except that instead of a ray in the
tangent space, the surfaces now have an embedded parameterised
disc at each marked point, together with a number $t \in [0,1/2]$.
We need to define on this space the structures above.  The circle
actions are given by rotating the discs. We need to say how to
glue two marked points together. If these have numbers $t,t' \in
[0,1/2]$ where $0 < t \le t'$, we glue together the circles of
radius $t$ around the marked points. If $t = 0$, we glue the
marked points together to get a node, with a ray in the tensor
product of the tangent lines for each side.

We get a chain of weak equivalences as follows. Let $\mc P_t$ be the
part  of $\mc P$ where all marked points have the same label $t \in
[0,1/2]$. The inclusion $\mc P_t \into \mc P$ is a weak equivalence.
Also there is a weak equivalence $\mc P_0 \to \mc N$, and a weak
equivalence $\mc M^s \to \mc P_{1/2}$.

\end{proof}

\begin{lemma}
The associated maps $\Fk(\mc M^s) \to \Fk(\mc P) \from \Fk(\mc N)$
of BV algebras induce isomorphisms on $\dhat$ homology.
\end{lemma}
\begin{proof}
We know the maps induce isomorphisms on $\d$ homology. The operator
$\dhat$ respects the grading by Euler characteristic of each of the
BV algebras.  On each graded piece, there are bounded filtrations by
number of marked points and number of connected components.  A
spectral sequence argument allows us to conclude the result.
\end{proof}

Let
$$
\tilmod (n) \to \cmod (n)
$$
be the principal $(S^1)^n$ orbi-bundle of curves in $\cmod(n)$
together with at each marked point a ray in the tangent space. There
is a map
$$
\mc N(n) \to \tilmod(n)
$$
which intertwines the $S^1 \wr S_n$ action.  This induces a map
$$
C_\ast(\mc N(n))_{hS^1 \wr S_n} \to C_\ast(\tilmod(n))_{hS^1 \wr
S_n}
$$
Let us redefine $\Fk(\cmod)$ as
$$
\Fk(\cmod) = \oplus_n  C_\ast(\tilmod(n))_{hS^1 \wr S_n}
$$
Evidently this is quasi-isomorphic to the previous definition.

By definition,
$$
\Fk(\mc N) = \oplus_n  C_\ast(\mc N(n))_{hS^1 \wr S_n}
$$
so that we have an algebra homomorphism $\Fk(\mc N) \to
\Fk(\cmod)$.  It is clear that this intertwines the ordinary
differential $\d$ on $\Fk(\mc N)$ with that on $\Fk(\cmod)$.
\begin{lemma}
The map $\pi : \Fk(\mc N) \to \Fk(\cmod)$ intertwines the
quantised differential $\dhat$ on $\Fk(\mc N)$ with the usual
differential $\d$ on $\Fk(\cmod)$.
\end{lemma}
\begin{proof}
Recall we can write
$$
\dhat = \d +  \tr
$$
where $\tr$ is an order $2$ operator on $\Fk(\mc N)$. It suffices to
show that $\pi(\tr (x)) = 0$ for all $x \in \Fk(\mc N)$.  Recall the
explicit description of $\tr$ at the end of section \ref{section
weyl fock}.

There is a gluing map $G_{ij} : \tilmod(n) \to \tilmod(n-2)$, for
each $1 \le i < j \le n$.  The diagram
$$
\xymatrix{ \mc N(n) \ar[r] \ar[d]^{G_{ij}} & \tilmod(n)
\ar[d]^{G_{ij}} \\
\mc N(n-2) \ar[r] & \tilmod(n-2) }
$$
commutes.

Also, the diagram
$$
\xymatrix{ \tilmod (n) \times S^1 \ar[r]^{a_i} \ar[d]^{p} &
\tilmod(n)
\ar[d]^{G_{ij}} \\
\tilmod(n) \ar[r]^{G_{ij}} & \tilmod(n-2) }
$$
commutes, where $p$ is the projection map and $a_i$ is the action
map for the $S^1$ action on $\tilmod(n)$ which rotates the ray at
the $i$'th marked point.

To show $\pi(\tr(x)) = 0$, it suffices to show the following. Let
$\op{D}_i$ be the degree one operator on $C_\ast(\tilmod(n))$
coming from the $i$'th circle action.  We need to show that, for
all $y \in C_\ast(\tilmod(n))$, $G_{ij} ( \op D_i y ) = 0$. That
this is a sufficient condition follows from the explicit
description of $\tr$ at the end of section \ref{section weyl
fock}.  But this condition follows from the commutativity of the
diagram
$$
\xymatrix{  C_\ast(\tilmod (n)) \otimes C_\ast (S^1) \ar[r]^{a_i}
\ar[d]^{p} & C_\ast(\tilmod(n))
\ar[d]^{G_{ij}} \\
C_\ast(\tilmod(n)) \ar[r]^{G_{ij}} & C_\ast( \tilmod(n-2)) }
$$

\end{proof}

As the BV operator (the differential) on $\Fk(\cmod)$ is a
derivation, the Poisson bracket associated to the BV structure is
trivial.  Therefore, the image of $S$ in $\Fk(\cmod)$ is closed.
We can write
$$
S = \sum \lambda^{2g-2+n} S_{g,n} \in \Fk(\mc N)[[\lambda]]
$$
The image $\pi(S_{g,n})$ of each $S_{g,n}$ is closed.

\begin{theorem}
The class $[\pi(S_{g,n})] \in H_{6g-6+2n}(\cmod(n))$ is the orbifold
fundamental class $[\cmod_{g,n}/S_n]$. \label{theorem dm master}
\end{theorem}
\begin{proof}
Filter $\Fk(\mc N)$ by saying $F^k (\Fk(\mc N))$ is the subspace
spanned  by chains on the space of surfaces with at least $k$ nodes.
Then the operator $\tr$  is of degree $1$ with respect to this, that
is $\tr(F^k) \subset F^{k+1}$.  It follows that $\{F^k,F^l\} \subset
F^{k+l+1}$.

We can similarly filter $\Fk(\cmod)$.   Note that
$$
H_\ast( \op{Gr}^0 \Fk (\cmod)  )  = H_\ast(\op{Gr}^0 \Fk(\mc N)) =
\oplus_n H_\ast( \cmod(n)/S_n,\partial \cmod(n)/S_n  )
$$
where $\partial$ refers to the boundary of the moduli space, i.e.\
the locus of nodal curves.

Let us use the notation
$$X(n) = \mc N(n) / S^1 \wr S_n$$
This is the moduli space of algebraic curves $C \in \cmod(n)/S_n$,
together with at each node, a ray in the tensor product of the
tangent spaces at each side.  Let $X_g(n) \subset X(n)$ be the
subspace of connected genus $g$ curves.

Let $\partial_i X(n) \subset X(n)$ be the locus with at least $i$
nodes. Then,
$$
H_\ast ( \op{Gr}^i \Fk(\mc N)) ) = \oplus_n H_\ast ( \partial_i
X(n),\partial_{i+1} X(n) )
$$

Since $S = \sum \lambda^{2g-2+n} S_{g,n}$ satisfies the master
equation, the classes $S_{g,n}$ are closed in $\op{Gr}^0 \Fk(\mc
N)$.  It suffices to show that the class
$$
[S_{g,n}] \in H_{6g-6+2n} ( X_g(n),\partial X_g(n))
$$
is the fundamental class of the orbifold with boundary $X_g(n)$.

The operator $\tr$, on homology, gives a map
$$
 \tr : H_\ast(\op{Gr}^0 \Fk(\mc N)) \to H_{\ast+1}(\op{Gr}^1
\Fk(\mc N))
$$
This corresponds to a map
$$\tr : H_\ast(X(n), \partial X(n) ) \to H_{\ast+1}( \partial X(n-2),\partial_2
X(n-2))$$ This  map has the following geometric description.  The
space $\partial X(n) \setminus \partial_2 X(n)$ is the moduli space
of curves $C \in X(n)$ with a single node.  Equivalently, it is the
moduli space of curves $C \in X(n+2)$ with a choice of two unordered
points, and a ray in the tensor product of the tangent spaces of
these points.  Thus there is a map
$$
\phi : \partial X(n) \setminus \partial_2 X(n) \to X(n+2) \setminus
\partial X(n+2)
$$

The space $X(n)$ has a canonical $\Q$ orientation, as it is a
complex orbifold.  This induces an orientation on the boundary
$\partial X(n)$. We can use Poincar\'e duality to identify
\begin{align*}
H_\ast( X(n),
\partial(X(n)) &= H^\ast(X(n) \setminus \partial
(X(n)) \\
H_\ast(\partial X(n),
\partial_{2}(X(n)) &= H^\ast(\partial X(n) \setminus \partial_{2} (X(n))
\end{align*}
Then the required map is minus the pull back map:
$$
\tr = -\phi^\ast : H^\ast  (  X(n+2) \setminus \partial (X(n+2) )
\to H^\ast (
\partial X(n) \setminus \partial_2 (X(n) )
$$

There is a natural boundary map $\d : H_\ast( X(n),\partial X(n))
\to H_{\ast-1}(\partial X(n),\partial_2 X(n))$.  The operator $\d +
 \tr$ makes $H_\ast (\Gr^{\le 1} \Fk(\mc N))$ into a BV algebra.
The algebra structure  is given by identifying it with $H_\ast(\Gr
\Fk(\mc N) / \Gr^{\ge 2} \Fk (\mc N))$.

This BV algebra structure gives us a Lie bracket
$$
\{\quad \}   : H_\ast(X(n),\partial X(n)) \otimes
H_\ast(X(m),\partial X(m)) \to H_\ast(\partial X(n+m-2),
\partial_2 X(n+m-2))
$$
This bracket has a similar geometric picture to the operator $\tr$.
We have a connected component $U \subset \partial X(n+m-2) \setminus
\partial_2 X(n+m-2)$ of curves where the node separates into two
components, one with $n$ and one with $m$ marked points. Then there
is a map $U \to X(n) \times X(m)$, and using Poincar\'e duality as
before, we can identify $\{\quad\}$ as minus the pull back in
cohomology.

The series $[S] = \sum \lambda^{2g-2+n} [S_{g,n}]$ satisfies the
master equation. This means that
$$
\d [S_{g,n}] + \tr [S_{g-1,n+2}] + \frac{1}{2} \sum_{\substack
{g_1+g_2 = g\\  n_1+n_2 = n+2}} \{ [S_{g_1,n_1}], [S_{g_2,n_2}] \}
=0
$$
It is clear that $[S_{g,n}] \in H_{6g-6+2n} (\cmod_{g,n}/S_n,
\partial \cmod_{g,n}/S_n )$ is the unique solution to this
equation, with the initial condition that $S_{0,3}$ is the
fundamental class, i.e. $\tfrac{1}{6}$ times the class of a point.

On the other hand, let $[X(n)] \in H_*(X(n),\partial X(n))$ be the
fundamental class.  It is clear that
$$\d [X(n)] + \tr [X(n+2)] = 0$$
Indeed, by the definition of the orientation on $\partial X(n)$,
\begin{align*}
\d [X(n)] &= [\partial X(n)] \in H_*(\partial X(n), \partial_2 X(n))  \\
\tr [X(n+2)] &= - [\partial X(n)] \in H_*(\partial X(n), \partial_2
X(n))
\end{align*}
Let $[X_g(n)]$ be the fundamental class of $X_g(n)$.  We have (in
some completion)
$$
\sum [X(n)] = \exp (\sum [X_g(n)])
$$
The equation $\d[X(n)] + \tr  [X(n+2)] = 0$ implies that
$$
\d [X_g(n)] +  \tr [X_{g-1,n+2}] + \frac{1}{2} \sum_{\substack
{g_1+g_2 = g\\  n_1+n_2 = n+2}} \{ [X_{g_1,n_1}], [X_{g_2,n_2}] \}
=0
$$
It is clear that $[X_{0,3}] = [S_{0,3}]$; in fact we chose $S_{0,3}$
to satisfy this.  It follows that $[X_{g,n}] = [S_{g,n}]$ for all
$g$ and $n$.
\end{proof}

Note that the space $C_\ast(\tilmod(n))_{h(S^1)^n}$ has an action of
$\K[t_1,\ldots,t_n]$.  On homology,
$H_\ast(C_\ast(\tilmod(n))_{h(S^1)^n}) = H_\ast(\cmod(n))$. The
action of $\K[t_1,\ldots,t_n]$ is by cap product with minus the
$\psi$ classes, i.e.\  the first Chern classes of the tangent lines
at the marked points. This is because the oriented torus bundle
$\tilmod(n) \to \cmod(n)$ is that associated to the tautological
tangent line bundles.

The space $C_\ast(\mc N)_{h(S^1)^n}$ also carries an action of
$\K[t_1,\ldots,t_n]$, and the map $C_\ast(\mc N)_{h(S^1)^n} \to
C_\ast(\tilmod(n))_{h(S^1)^n}$ is a map of $\K[t_1,\ldots,t_n]$
modules.  Therefore we can see not just the fundamental class, but
its cap products with $\psi$ classes.

\section{The Gromov-Witten type invariants associated to a TCFT}
\label{section gw invariants} In this section, for each TCFT we will
construct a left ideal in the associated Weyl algebra, together with
an element in it. This encodes the Gromov-Witten potential  of the
TCFT.

Recall that so far, for each TCFT $F$, we have constructed a Weyl
algebra $\mc W(F)$, and a Fock space $\Fk (F)$.  There is also Weyl
algebra and Fock space, $\mc W(\mc M)$ and $\Fk(\mc M)$, associated
to moduli space.

As I mentioned earlier, in an ideal world, we might hope that there
is a natural transformation $C'_\ast(\mc M) \to F$ of functors from
$H_*(A) \to \op{Comp}_{\K}^{\Z/2}$, encoding the TCFT structure on
$F$.  However this is not the case, as  only moduli spaces of
surfaces with at least one incoming boundary act.

However, suppose there was such a natural transformation
$C'_\ast(\mc M) \to F$.  Then there would be an associated algebra
homomorphism $\mc W(\mc M) \to \mc W(F)$, and a map $\Fk(\mc M) \to
\Fk(F)$ of $\mc W(\mc M)$ modules.  The closed element $\exp(S) \in
\Fk(\mc M)[[\lambda]]$ then maps to a closed element in the Fock
space $\Fk(F)[[\lambda]]$.

We want to mimic this construction in the real-world situation.
Recall that then, instead of a natural transformation $C'_\ast(\mc
M) \to F$, we have a twisted functor $F^{\mc M}$, with
transformations $C'_\ast(\mc M) \to F^{\mc M} \from F$.   Then
$\Fk(F^{\mc M})$ is a $\mc W(\mc M)-\mc W(F)$ bimodule, and there is
a map $\Fk( \mc M ) \to \Fk(F^{\mc M})$ of $\mc W(M)$ modules.

The closed element $\exp(S) \in \Fk (\mc M)[[\lambda]]$ now passes
to a closed element in $\Fk(F^{\mc M})[[\lambda]]$, which we call
$\mc D$. That is, we have constructed a module for the Weyl algebra
$\mc W(F)[[\lambda]]$ and an element in it.

Of course, the space $\Fk(F^{\mc M})$ is far too big.  As an algebra
it is isomorphic to $\Fk(\mc M) \otimes \Fk(F)$.  The differential
does not respect this decomposition, nor does the action of $\mc
W(F)$.

Let
$$\mc W_-(F)  = \oplus_n F_{hS^1}(n)_{S_n} \subset \mc W(F)$$ be the commutative subalgebra of
elements which contain only negative powers of the $t_i$. This
subalgebra is not preserved by the differential.

\begin{lemma}
The action of $\mc W_-(F)[[\lambda]]$ on $\mc D \in \Fk(F^{\mc
M})[[\lambda]]$ generates a free $\mc W_-(F)[[\lambda]]$ submodule,
which we call $\Fk_{\mc D}(F)$. Further, $\Fk_{\mc D}(F)$ is also
preserved by the action of all of $\mc W(F)[[\lambda]]$, and by the
differential. \label{lemma submodule}
\end{lemma}
\begin{proof}
This is a matter of unravelling the definitions of the various
algebras and their modules.  The element $\mc D \in \Fk(F^{\mc
M})[[\lambda]]$ is of the form $\exp(S)$, where
$$S \in \left( \oplus \Fk^{g,n}(\mc M)\right) [[\lambda]] \subset \Fk(F^{\mc
M})[[\lambda]]$$ Then, the submodule $\Fk_{\mc D}(F)$ is explicitly
given as the set
$$
\exp(S) \otimes f
$$
where $f \in \Fk(F)[[\lambda]]$.  The only thing to check is that
this space is preserved by the action of $\mc W(F)$.    If we have
 an element  $X \in F(1)[[t]] \subset \mc W(F)$, then it acts by a
derivation of $\Fk(F^{\mc M})$. This derivation takes elements of
$\Fk^{g,n}(\mc M)$ into $\Fk(F)$.  Therefore
$$
X( \exp(S) \otimes f ) = \exp(S) \otimes X f + \exp(S) \otimes (XS)
f
$$
If the map $F(1)^{\otimes n} \to F(n)$ was an isomorphism, and not
just a quasi-isomorphism, this would be enough.  In the general
case, a little bit more needs to be checked, but its not difficult.
\end{proof}

By definition, the action of $\mc W(F)[[\lambda]]$ on $\Fk_{\mc
D}(F)$ generates the module from the element $\mc D$.  The
annihilator of $\mc D$ is a left ideal in $\mc W(F)[[\lambda]]$,
which should be viewed as the fundamental object encoding the
Gromov-Witten type potential associated to a TCFT.

\section{Choice of polarisation and Givental's group}

To extract a more familiar looking potential, we need to pass to
homology. In this section we will make the following two assumptions
on our TCFT:
\begin{enumerate}
\item
The map $\op D : F(1) \to F(1)$ is zero on homology.  When our TCFT
comes from a Calabi-Yau $A_\infty$ category as in \cite{cos_2004oc},
this is equivalent to the degeneration of the spectral sequence from
Hochschild homology tensored with $\K((t))$ to periodic cyclic
homology. This spectral sequence is the non-commutative analog of
Hodge to de Rham.
\item
The inner product on $H_\ast(F(1)))$ is non-degenerate.
\end{enumerate}
This implies that $H_\ast(F(1)_\Tate)$ is symplectic, i.e.\ the
natural anti-symmetric pairing is non-degenerate.  Also the map
$H_\ast(F(1)^{hS^1}) \to H_\ast(F(1)_\Tate)$ is injective, and the
image is a Lagrangian subspace.

In this section we will use the following notation :
\begin{align*}
\mc H &= H_\ast(F(1)_\Tate)\\
H &= H_\ast(F(1)) \\
\mc H_+ &= H_\ast(F(1)^{hS^1}) \subset \mc H
\end{align*}
We have the Weyl algebra $\mc W(\mc H)$ and the Fock space $\Fk(\mc
H)$, which is the quotient of $\mc W(\mc H)$ by the left ideal
generated by $\mc H_+$.

\begin{lemma}
There  are natural isomorphisms
\begin{align*}
H_\ast(\mc W(F)) &\iso \mc W(\mc H) \\
H_\ast(\Fk(F)) & \iso \Fk (\mc H )
\end{align*}
The second isomorphism is compatible with the $\mc W(\mc H)$
actions.
\end{lemma}
\begin{proposition}
There is an isomorphism of $\mc W(\mc H)[[\lambda]]$ modules
$$
\mf F(\mc H)[[\lambda]] \iso H_\ast(\mf F_{\mc D}(F))
$$
which is unique up to multiplication by an  element of
$\K[[\lambda]]^\times$.
\end{proposition}
\begin{proof}

The uniqueness part is well known.  We will show  existence.

Modulo $\lambda$, this follows immediately from the definition of
$\mf F_{\mc D}(F)$. The point is that the potential $\mc D =
\exp(S)$ is $1$ modulo $\lambda$.

The $\mc W(H)$  module $\mf F(\mc H)$ is rigid, meaning that any
flat deformation of it over a pro-nilpotent local ring is trivial.
It remains to show that $H_\ast(\mf F_{\mc D}(F))$ is flat over
$\K[[\lambda]]$.

Let $\mc H_- \subset \mc H$ be any complementary subspace to $\mc
H_+$, so that $\mc H = \mc H_- \oplus \mc H_+$.  The action of $\mc
W(\mc H)[[\lambda]]$ on $\mc D \in H_\ast(\mf F_{\mc D}(F))$ gives a
map
$$
\Sym^\ast \mc H_- [[\lambda]] \to H_\ast( \mf F_{\mc D}(F))
$$
We need to show that this is an isomorphism. Indeed, since  the
operator $\op D : F(1) \to F(1)$ is zero on homology,  we can lift
the subspace $\mc H_-$ to an isotropic subcomplex $\mc H_- \subset
F(1)_{hS^1}\subset F(1)_\Tate$.  This makes $\Sym^\ast \mc H_-$  a
subalgebra of $\mc W(\mc H)$, preserved by the differential. The map
$\Sym^\ast \mc H_-[[\lambda]] \to \Fk_{\mc D}(F)$ given by acting on
$\mc D \in \Fk_{\mc D}(F)$ is compatible with the differentials. The
differential on $\Fk_{\mc D}(F)$ can be written as $\d + \tr$, where
$\d$ is the usual differential on $\oplus
F_{hS^1}(n)_{S_n}[[\lambda]]$ and $\tr$ is an order $2$ operator.
The operator $\tr$ is not the usual BV operator, but incorporates
the action of $\mc M$ on $F$ and the solution of the master equation
$S \in \Fk(\mc M)[[\lambda]]$. However, $\tr$ contains the circle
operator $\op D$.  Therefore, on $\d$ homology, $\tr$ is zero.

We know from lemma \ref{lemma submodule} that the map $\Sym^\ast \mc
H_-[[\lambda]] \to \Fk_{\mc D}(F)$ induces an isomorphism on $\d$
homology.  This implies it induces an isomorphism on $\d + \tr$
homology.
\end{proof}

This proposition shows that we have a canonically defined line (i.e.
rank one $\K[[\lambda]]$ submodule)\footnote{We can use the dilaton
equation to reduce the ambiguity from $\K[[\lambda]]^{\times}$ to
$\K^\times$.} $\ip{\mc D} \subset \mf F[[\lambda]]$. This plays the
role of the total ancestral potential.

In order to get a more familiar kind of potential, we need to choose
some extra data.  The symplectic vector space $\mc H$ has various
natural structures, namely the isomorphisms given by multiplication
by $t$ and $t^{-1}$, the Lagrangian subspace $\mc H_+$.  We will
look for polarisations of $\mc H$ compatible with these structures.
\begin{definition}
A compatible polarisation of $\mc H$ is a Lagrangian subspace $\mc
H_- \subset \mc H$ such that
\begin{enumerate}
\item
$$
\mc H = \mc H_- \oplus \mc H_+
$$
\item
The operator $t^{-1}$ preserves $\mc H_-$.
\end{enumerate}
\end{definition}
The space $\mc H$ is naturally filtered, by the subspaces $t^k \mc
H_+$.  The associated graded is canonically isomorphic to $H((t))$,
as a $\K((t))$ module.  The corresponding symplectic form on
$H((t))$ is given by
$$
\Omega(a \otimes f,b \otimes g) = \ip{a,b}\op{Res} f(-t)g(t) \op d t
$$
where $\ip{\quad}$ refers to the pairing on $H = H_*(V)$ coming from
that on $V$.

\begin{lemma}
A compatible polarisation of $\mc H$ is the same data as an an
isomorphism of symplectic vector spaces
$$
\mc H \iso H((t))
$$
which is compatible with the action of $t$, takes $\mc H_+$ to
$H[[t]]$, and on the associated graded spaces is the identity.
\end{lemma}
\begin{proof}
The spaces $t^k \mc H_-$ give a splitting of the filtration.
\end{proof}

If we pick a compatible polarisation, then we can identify the Fock
space $\mf F$ with $\Sym^\ast \mc H_-$. $\mc W(\mc H)$ acts on
$\Sym^\ast \mc H_- $ in a standard fashion; vectors in $\mc H_-
\subset \mc H$ act by multiplication, while those in $\mc H_+
\subset \mc H$ act by differentiation by the corresponding element
of $\mc H_-^{\ast}$.

Since $\mc H_- = t^{-1} H[t^{-1}]$,  for each choice of compatible
polarisation of $\mc H$, we find a line
$$
\ip{\mc D} \subset \left( \Sym^\ast t^{-1} H[t^{-1}]
\right)[[\lambda]]
$$

Suppose we change the polarisation of $\mc H$.  This corresponds to
a change of the isomorphism
$$
\mc H \iso H((t))
$$
Any such is given by a symplectomorphism of $H((t))$. When we change
this isomorphism, $\mc H_+$ must again correspond to $H[[t]]$, and
the new isomorphism must again be compatible with the action of $t$,
and be the identity on the associated graded. This implies that the
corresponding symplectomorphism of $H((t))$ is of the form
$$
\phi(t) =  1 + \sum_{t \ge 1} t^i \phi_i
$$
where $\phi_i : H \to H$ are linear maps.  The symplectomorphism
condition translates into
$$
\phi(t)\phi^\ast(-t) = 1
$$
where $\phi^\ast$ is obtained by replacing each $\phi_i$ by its
adjoint.  This means that $\phi(t)$ is an upper-triangular element
in Givental's twisted loop group. \footnote{If we allowed a change
of $\mc H_+$
 we would find an element of the full twisted loop group, not
just the upper triangular part. However, it is more difficult  to
make sense of this outside the formal neighbourhood of the
upper-triangular part of the twisted loop group.  }

Any such operator can be quantised  to act on the Fock space.  As,
any such $\phi$ admits a logarithm, which is an infinitesimal
symplectomorphism of $H((t))$, commuting with the $t$ action. Such
an infinitesimal symplectomorphism, $A$ say,  admits a quantisation,
in a standard fashion \cite{giv2001}, to an element in the Weyl
algebra $\what{A} \in \mc W(H((t)))$. $\what{A}$ is characterised up
to an additive scalar by the condition that the inner derivation
$[\what{A},-]$ of $\mc W(H((t)))$, when restricted to $H((t))
\subset \mc W(H((t)))$, is $A$.

Then the quantised operator $\what{\phi}$ acts on the Fock space by
$$
\what{\phi} = \exp( \what{\log(\phi) } )
$$
where $\what{\log(\phi)}$ acts on the Fock space $\Sym^{\ast} t^{-1}
H[t^{-1}]$ in the standard way, as an element of the Weyl algebra.
The exponential makes sense because $\what{\log(\phi)}$ is a locally
nilpotent operator on $\Sym^\ast t^{-1} H[t^{-1}]$.

The symplectomorphism $\phi$ induces an automorphism $\phi$ of the
Weyl algebra $\mc W(H((t)))$.  We can twist the action of $\mc
W(H((t)))$ on $\Sym^\ast t^{-1} H[t^{-1}]$ via the automorphism
$\phi$.
\begin{lemma}
$\what{\phi}$ is the unique up to scale isomorphism of $\Sym^\ast
t^{-1} H[t^{-1}]$ such that for $w \in \mc W(\mc H((t)))$, $x \in
\Sym^\ast t^{-1} H[t^{-1}]$,
$$
\what{\phi} ( w \cdot x) = \phi (w)\cdot \what{\phi} (x)
$$
Therefore, if we change the isomorphism $\mc H \iso H((t))$ by a
symplectomorphism $\phi$, then the corresponding line $\ip{\mc D}$
changes by $\what{\phi}$.
\end{lemma}
\begin{proof}
Let $A = \log \phi$. Then
$$
\phi(w) = \exp ( \op {ad} \what A )  ( w)
$$
and
\begin{align*}
\what{\phi} ( w \cdot x ) &= \exp ( \what A ) \cdot w \cdot x \\
&= \exp (\what A) \cdot w \cdot \exp( -\what A ) \cdot \exp(\what A) \cdot x \\
&= \phi(w) \cdot \what{\phi} (x)
\end{align*}

\end{proof}

\section{Relation with ordinary Gromov-Witten invariants}
Suppose our TCFT is equipped with operations from the
Deligne-Mumford spaces, and not just the uncompactified spaces. This
should happen for the $A$-model TCFT associated to a compact
symplectic manifold. Then the ancestral potential constructed above,
for a certain choice of polarisation, coincides with the actual
ancestral potential, coming from the fundamental class and $\psi$
classes in Deligne-Mumford space.

In fact, this is immediate from the results of section \ref{section
fundamental master}, but I will briefly go through some of the
details anyway.

For simplicity we will suppose that the TCFT is also equipped with
operations from the space of curves with no marked points, all
though this is not necessary.  Also we will assume, to keep the
notation simple, that our TCFT is split.

So, suppose we are given a complex $V$, with an inner product, and
maps
$$
\phi : C_\ast(\cmod(n)) \to V^{\otimes n}
$$
such that the gluing maps between the spaces $\cmod(n)$ correspond
to the inner product on $V$.

Give $V$ the trivial circle action, $\op D = 0$.  Then $V$ defines a
functor, say  $F: H_\ast(A) \to \op{Comp}_\K$.

Recall the space $\tilmod(n)$ is the principal $(S^1)^n$ orbi-bundle
over $\cmod(n)$ given by curves in $\cmod(n)$ with at each marked
point a ray in the tangent space. The complexes $C_\ast(\tilmod(n))$
define a functor $H_\ast(A) \to \op{Comp}_\K$.  The complexes
$C_\ast(\mc N(n))$ do also, and there is a natural transformation
$$
C_\ast(\mc N) \to C_\ast(\tilmod)
$$
There is a natural transformation $C_\ast(\tilmod) \to F$.

There is an associated map
$$
\Fk(\mc N) \to \Fk(\tilmod) \to \Fk ( F  )
$$
We have
$$
\Fk( \tilmod) \simeq \oplus_n C_\ast(\cmod(n))_{S_n}
$$
and we have already seen that the image of the solution of the
master equation in $\Fk(\mc N)$ goes to the fundamental classes of
$\cmod_{g,n}/S_n$.

Now, since the circle operator on $V$ is zero, there is a canonical
isomorphism
$$
\Fk(F) =  \Sym^\ast t^{-1} V[t^{-1}]
$$
This arises from the canonical polarisation
$$V((t)) = V[[t]] \oplus t^{-1} V[t^{-1}]$$
which in this case is compatible with the differential.

The map
$$
C_\ast(\cmod(n))_{S_n} \to \Sym^n t^{-1} V[t^{-1}]
$$
is given, at least on homology, by
$$
x \mapsto \sum_{l_1,\ldots,l_n \ge 0} (-1)^{\sum l_i} t_1^{-l_1-1}
\ldots t_n^{-l_n-1} \phi(\psi_1^{l_1} \ldots \psi_n^{l_n}\cap x)
$$
The point is that the action of $\K[t_1,\ldots,t_n]$ on
$C_\ast(\cmod(n))\simeq C_\ast(\tilmod(n))_{h(S^1)^n}$ is given on
homology by cap product with minus the $\psi$ classes.

This makes it clear that the ray in the Fock space we have
constructed coincides with the ordinary ancestral potential.

\section{The B model and the holomorphic anomaly}
\label{section B model}

Suppose that $X$ is a smooth projective Calabi-Yau variety of
dimension $d$, for simplicity  over  $\C$ (which in this section we
take to be our base field). Pick a holomorphic volume form
$\op{Vol}_X$ on $X$. We would like to use the results of this paper
and \cite{cos_2004oc} to construct the B model analog of the
Gromov-Witten potential of $X$. Currently there is a small technical
gap in this construction, which will be discussed elsewhere.

In \cite{cos_2004oc}, I showed how to associated to a Calabi-Yau
$A_\infty$ category a TCFT, whose homology is the Hochschild
homology of the category.  The derived category of coherent sheaves
on $X$, $\mc D^b(X)$, is a Calabi-Yau category.  However, it is
\emph{not} the category we want, for various reasons explained in
\cite{cos_2004oc}.

Something a bit closer to what we want is a differential graded
category of complexes of sheaves on $X$.  There are various
versions, all of which should be quasi-isomorphic.  Perhaps the
simplest is to consider the category whose objects are bounded
complexes of (algebraic) vector bundles on $X$, and whose complex of
morphisms $E \to F$ is the Dolbeaut resolution of $E^\ast \otimes
F$.  Other constructions, e.g.\ using Cech resolutions, or any
injective resolution functor, have the advantage of working over
fields other than $\C$.

This dg model for the category of sheaves on $X$ does not quite
satisfy the conditions of \cite{cos_2004oc} either.  The homology of
this dg category is the derived category.  What we need to do is to
give the derived category an $A_\infty$ structure, using the
categorical analog of Kadeishvili's theorem \cite{kad1982}. In order
for the resulting $A_\infty$ category to be of Calabi-Yau type, the
higher products need to be cyclically symmetric with respect to the
pairing.

As far as I am aware, no-one has written down a proof that the
higher products can be made cyclically symmetric, and that the
resulting $A_\infty$ category has the ``correct'' Hochschild
homology.  This is the technical gap mentioned earlier.  Probably
one could use Hodge theory, and the explicit form of the homological
perturbation lemma \cite{kon_soi2000}, to prove the first part.

Let us suppose that one can do this, and denote by $\mc
D_\infty^b(X)$ the resulting Calabi-Yau $A_\infty$ category.   To
this, the results of \cite{cos_2004oc} associate a TCFT, $F$.  The
homology of this is the Hochschild homology of $\mc D_\infty^b(X)$.
We have:
\begin{align*}
H_i(F(1)) &= HH_{i-d}(\mc D_{\infty}^b(X))\\
H_i(F(1)_{hS^1}) &= HC_{i-d}(\mc D_\infty^b(X))\\
H_i(F(1)^{hS^1}) &= HC_{i-d}^{-} (\mc D_\infty^b(X)) \\
H_i(F(1)_\Tate) &= HP_{i-d}(\mc D_\infty^b(X))
\end{align*}
Here, $HH_\ast$ is Hochschild homology, $HC_\ast$ is cyclic
homology, $HC_\ast^{-}$ is negative cyclic homology, and $HP_\ast$
is periodic cyclic homology. The shift in grading is due to a
difference in grading conventions between this paper and
\cite{cos_2004oc}.

All of these groups can be identified with more classical cohomology
groups of $X$.  Let $H^{-\ast}(X,\C)$ be the usual cohomology of
$X$, with the grading reversed.    Let $F^p H^{-\ast}(X,\C)$ be the
$p$'th part of the Hodge filtration, i.e. the part coming from
$(r,s)$ forms where $r \ge p$.  We should have
\begin{align*}
HH_i(\mc D^b_\infty(X)) &= \oplus_{q - p = i} H^{p}(X,\Omega_X^q) \\
HP_\ast(\mc D^b_\infty(X)) &= H^{-\ast}(X,\C) \otimes \C((t)) \text{
where } t \text{ has degree } -2\\
HC_\ast^{-}(\mc D^b_\infty(X)) &= \oplus_p  F^p H^{-\ast}(X,\C)
\otimes t^{-p} \C[[t]]\\
HC_\ast(\mc D^b_\infty(X)) &= HP_\ast(\mc
D^b_\infty(X))/HC_\ast^{-}(\mc D^b_\infty(X))
\end{align*}
To see this, consider the
$$ V_i = \oplus_{p-q=i -d } \Omega_X^{p,q}$$ with differential $\delbar$.     The operator
$\partial$ is a circle action. We should (if $\mc D_\infty^b(X)$ has
been constructed correctly) have a HKR quasi-isomorphism
$$
V \simeq F(1)
$$
in such a way that the circle action $\partial$ on $V$ corresponds
to $\op D$ on $F(1)$.   Also, $F(1)$ is quasi-isomorphic to the
 Hochschild chain complex shifted by $d$, where the circle operator
corresponds to the B operator. Then  $ F(1)_\Tate \simeq V_\Tate $,
$V_{hS^1} \simeq F(1)_{hS^1}$ and $V^{hS^1} \simeq F(1)^{hS^1}$.

Now consider the de Rham complex $\Omega_X^{-\ast,-\ast}$, with the
reverse grading and the usual differential $\op d = \partial +
\delbar$. Define a map
$$
\Omega_X^{-\ast,-\ast} ((t)) \to V_\Tate
$$
which is $\C((t))$ linear and sends $\alpha \in \Omega_X^{p,q}$ to
$t^p \alpha$, using the identification $V_i = \oplus_{p-q=i-d}
\Omega_X^{p,q}$.

This is clearly an isomorphism of complexes, which shifts degree by
$d$. So we find that
$$
H_{\ast + d}(V_\Tate) = H^{-\ast}(X,\C)((t))
$$
The subcomplex $V^{hS^1}$ corresponds to the subcomplex of
$\Omega_X^{-\ast,-\ast}((t))$ spanned by $\alpha t^k$, where $\alpha
\in \Omega_X^{p,q}$ and $k+p \ge 0$.  Thus, $H_\ast(V^{hS^1})$
corresponds to $\oplus F^p H^{-\ast}(X,\C) \otimes t^{-p} \C[[t]]$.

The material in this paper shows  (with the caveat discussed
earlier) that there is a canonically defined line in a Fock space
for this vector space, with a certain symplectic form, which plays
the role of the Gromov-Witten potential.

\subsection{The holomorphic anomaly}

In \cite{bcov}, Bershadsky et al. show that the B model potential
does not vary holomorphically on moduli space, but has an
``anomaly''. This was reinterpreted by Witten \cite{wit1993} to say
that the B model potential is an element of the Fock space modelled
on the symplectic vector space $H^3(X)$.  This vector space has a
Gauss-Manin flat connection on the Calabi-Yau moduli space, and the
vector in the Fock space is flat for the associated projectively
flat connection.  The results of this paper fit in very well with
Witten's viewpoint on the holomorphic anomaly.

Let $\mc M_{CY}$ be moduli space of Calabi-Yau A-infinity categories
(whatever that means).  On $\mc M_{CY}$, we have a sheaf of Weyl
algebras, with a left ideal in it encoding the potential.  The sheaf
of Weyl algebras should have a natural connection, flat up to a
coherent system of homotopies. On homology this will be the
connection induced by the Gauss-Manin connection on periodic cyclic
homology.
\begin{conjecture}
After taking account of the unstable moduli spaces $(g,n) =
(0,1),(0,2)$, the ideal is preserved up to homotopy by the flat
connection.
\end{conjecture}
This will be discussed elsewhere.

Recall that in order to get a more familiar kind of potential, we
need to pick a polarisation of the symplectic vector space $\mc H$.
As always, one half of the polarisation is defined; we have the
subspace $\mc H_+$. Using the identification $\mc H =
H^{-\ast}(X)((t))$ from the previous subsection, we have seen that
$\mc H_+$ is defined using the Hodge filtration.  A complementary
subspace $\mc H_-$ can be constructed using a splitting of the Hodge
filtration.  There is of course a natural splitting, given by the
complex conjugate filtration.

Therefore, for each Calabi-Yau $X$, with choice of holomorphic
volume form, we get a $\C[[\lambda]]$ line
$$
\ip{\mc D_X } \subset (\Sym^\ast t^{-1} H^\ast(X) [t^{-1}])
[[\lambda]]
$$
Rescaling the holomorphic volume form corresponds to rescaling
$\lambda$.  The fact that the potential is flat tells us what
happens when we change the complex structure on $X$.  This
corresponds to  keeping the same symplectic vector space, Fock
space, and line in the Fock space, but changing the polarisation of
$\mc H$. 

\def\cprime{$'$}

\end{document}